\documentclass[a4paper]{article}

\usepackage{amsmath, amssymb, amsfonts, amsthm, bm}

\newcommand{\vol}[1]{\operatorname{vol}(#1)}
\DeclareMathOperator{\CMdet}{CMdet}
\DeclareMathOperator{\prob}{Prob}
\DeclareMathOperator{\rank}{rank}
\DeclareMathOperator{\trace}{trace}
\DeclareMathOperator{\sgn}{sgn}

\DeclareMathOperator{\conv}{conv}

\DeclareMathOperator{\interior}{int}
\DeclareMathOperator{\Span}{span}
\DeclareMathOperator{\binlog}{lg}
\newcommand{\half}{\frac12}

\newcommand{\epsi}{\varepsilon}
\newcommand{\fhi}{\varphi}
\newcommand{\Tr}{\mathrm{tr}}
\newcommand{\norm}[1]{\lVert#1\rVert}
\newcommand{\ipr}[2]{\left\langle #1, #2 \right\rangle}

\newcommand{\Sphere}{\mathbb{S}}
\newcommand{\R}{\mathbb{R}}

\newcommand{\N}{\mathbb{N}}

\newcommand{\vect}[1]{\bm{#1}}
\newcommand{\bvect}[1]{\overline{\bm{#1}}}
\newcommand{\hvect}[1]{\widehat{\bm{#1}}}
\newcommand{\vvect}[1]{\widetilde{\bm{#1}}}
\newcommand{\va}{\vect{a}}
\newcommand{\bva}{\bvect{a}}
\newcommand{\hva}{\hvect{a}}
\newcommand{\vva}{\vvect{a}}
\newcommand{\vb}{\vect{b}}
\newcommand{\bvb}{\bvect{b}}
\newcommand{\hvb}{\hvect{b}}
\newcommand{\vvb}{\vvect{b}}
\newcommand{\vc}{\vect{c}}

\newcommand{\ve}{\vect{e}}
\newcommand{\vf}{\vect{f}}
\newcommand{\vx}{\vect{x}}
\newcommand{\bvx}{\bvect{x}}
\newcommand{\hvx}{\hvect{x}}
\newcommand{\vvx}{\vvect{x}}
\newcommand{\vy}{\vect{y}}

\newcommand{\vs}{\vect{s}}
\newcommand{\vo}{\vect{o}}
\newcommand{\vp}{\vect{p}}
\newcommand{\bvp}{\bvect{p}}
\newcommand{\vq}{\vect{q}}

\newcommand{\vu}{\vect{u}}
\newcommand{\vv}{\vect{v}}
\newcommand{\vw}{\vect{w}}
\newcommand{\collection}[1]{\mathcal{#1}}

\newcommand{\CP}{\collection{P}}

\newcommand{\abs}[1]{\lvert#1\rvert}
\newcommand{\bigabs}[1]{\left\lvert#1\right\rvert}
\newcommand{\card}[1]{\lvert#1\rvert}

\newcommand{\length}[1]{\operatorname{length}(#1)}
\newcommand{\hit}[1]{\operatorname{hit}(#1)}

\theoremstyle{plain}
\newtheorem{theorem}{Theorem}
\newtheorem{lemma}{Lemma}
\newtheorem{corollary}{Corollary}
\newtheorem{conjecture}{Conjecture}
\newtheorem{problem}{Problem}
\newtheorem{proposition}{Proposition}
\newtheorem*{Dvoretzky}{Dvoretzky's Theorem}
\newtheorem*{ranklemma}{Rank Lemma}

\newtheorem*{JacksonV}{Jackson V}

\begin{document}

\bibliographystyle{amsplain}

\title{Equilateral sets in finite-dimensional normed spaces}
\author{Konrad J. Swanepoel\thanks{Department of Mathematics, Applied Mathematics and Astronomy, University of South Africa, PO Box 392, Pretoria 0003, South Africa.
	E-mail: \texttt{swanekj@unisa.ac.za}.
This material is based upon work supported by the National Research Foundation under Grant number 2053752. These are lecture notes for a $4$-hour course given on April 20-23, 2004, at the Department of Mathematical Analysis, University of Seville. I thank the University of Seville for their kind invitation and generous hospitality.}}
\date{}
\maketitle

\begin{abstract}
This is an expository paper on the largest size of equilateral sets in finite-dimensional normed spaces.
\end{abstract}

\section*{Contents}
\noindent \textbf{1\; Introduction}

1.1\; Definitions

1.2\; Kusner's questions for the $\ell_p$ norms

\medskip
\noindent \textbf{2\; General upper bounds}

2.1\; The theorem of Petty and Soltan

2.2\; Strictly convex norms

\medskip
\noindent \textbf{3\; Simple lower bounds}

3.1\; General Minkowski spaces

3.2\; $\ell _p^n$ with $1<p<2$

\medskip
\noindent \textbf{4\; Cayley-Menger Theory}

4.1\; The Cayley-Menger determinant

4.2\; Embedding an almost equilateral $n$-simplex into $\ell _2^n$

\medskip
\noindent \textbf{5\; The Theorem of Brass and Dekster}

\medskip
\noindent \textbf{6\; The Linear Algebra Method}

6.1\; Linear independence

6.2\; Rank arguments

\medskip
\noindent \textbf{7\; Approximation Theory: Smyth's approach}

\medskip
\noindent \textbf{8\; The best known upper bound for $e(\ell _1^n)$}

\medskip
\noindent \textbf{9\; Final remarks}

9.1\; Infinite dimensions

9.2\; Generalizations

\section{Introduction}
In this paper I discuss some of the more important known results on equilateral sets in finite-dimensional normed spaces with an analytical flavour.
I have omitted a discussion of some of the low-dimensional results of a discrete geometric nature.
The omission of a topic has no bearing on its importance, only on my prejudices and ignorance.
The following topics are outside the scope of this paper: almost-equilateral sets, $k$-distance sets, and equilateral sets in infinite dimensional spaces.
However, see the last section for references to these topics.

I have included the proofs that I find fascinating.
The arguments are mostly very geometrical, using results from convex geometry and the local theory of Banach spaces, and also from algebraic topology.
In the case of the $\ell_p$-norm the known results use tools from linear algebra, probability theory, combinatorics and approximation theory.

Equilateral sets have been applied in differential geometry to find minimal surfaces in finite-di\-men\-sio\-nal normed spaces \cite{MR95i:58051, MR93h:53012}.

Throughout the paper I state problems.
Perhaps some of them are easy, but I do not know the answer to any of them.

\subsection{Definitions}
A subset $S$ of a metric space $(X,d)$ is \emph{$\lambda$-equilateral} ($\lambda>0$) if $d(x,y)=\lambda$ for all $x,y\in S, x\neq y$.
Let $e_\lambda(X)$ denote the largest size of a $\lambda$-equilateral set in $X$, if it exists, otherwise set $e_\lambda(X)=\infty$.
If $X$ is compact then it is easily seen that $e_\lambda(X)<\infty$ for each $\lambda>0$.

In a normed space $X$, if $S$ is $\lambda$-equilateral, then $\frac{\mu}{\lambda}S$ is $\mu$-equilateral for any $\mu>0$; thus the specific value of $\lambda$ does not matter, and we only write $e(X)$.
If $X$ is finite-dimensional, then $e(X)<\infty$.
In the remainder of the paper we abbreviate finite-dimensional normed space to \emph{Minkowski space}.
The \emph{Banach-Mazur distance} between two Banach spaces $X$ and $Y$ is the infimum of all $c\geq 1$ such that there exists a linear isomorphism $T:X\to Y$ satisfying $\norm{T}\norm{T^{-1}}\leq c$.
In the finite-dimensional case, this infimum is of course attained.
We always use $n$ to denote the dimension of $X$, and $m$ the size of a set.
We use $\R$ for the field of real numbers, $\R^n$ for the $n$-dimensional real vector space of $n$-dimensional column vectors, and $\ve_i$ for the $i$th \emph{standard basis vector} which has $1$ in position $i$ and $0$ in all other positions.
The transpose of a matrix $A$ is denoted by $A^\Tr$.
For a set $S\subseteq \R^n$, $\conv(S), \interior(S), \vol{S}$ denotes its convex hull, interior, and Lebesgue measure (if $S$ is measurable), respectively.
By $\log x$ we denote the natural logarithm, and by $\binlog x$ the logarithm with base $2$.

\subsection{Kusner's questions for the $\ell_p$ norms}\label{kusner}
On $\R^n$ the $\ell_p$ norm is defined by
\[ \norm{(x_1,x_2,\dots,x_n)}_p := \left(\sum_{i=1}^n\abs{x_i}^p\right)^{1/p}\]
if $1<p<\infty$, and for $p=\infty$ by
\[ \norm{(x_1,x_2,\dots,x_n)}_\infty := \max\{\abs{x_i} : i=1,\dots,n\}.\]
We denote the Minkowski space $(\R^n, \norm{\cdot}_p)$ by $\ell_p^n$.
The $\ell_p$ distance between two points $\va, \vb\in\R^n$ is $\norm{\va-\vb}_p$.
Thus Euclidean space is $\ell_2^n$.

The following lower bounds for $e(\ell_p^n)$ are simple:
\begin{itemize}
\item Since $\{\pm\ve_i:i=1,\dots,n\}$ is $2$-equilateral in $\ell_1^n$, we have $e(\ell_1^n)\geq 2n$.
\item For any $p\in(1,\infty)$ and an appropriate choice of $\lambda=\lambda(p)\in\R$, the set
\[\{\ve_1,\ve_2,\dots,\ve_n,\lambda\sum_{i=1}^n\ve_i\}\]
 is $2^{1/p}$-equilateral in $\ell_p^n$; hence $e(\ell_p^n)\geq n+1$.
\item $\{(\epsi_1,\dots,\epsi_n):\epsi_i=\pm 1\}$ is $2$-equilateral in $\ell_\infty^n$; thus $e(\ell_\infty^n)\geq 2^n$.
\end{itemize}

A theorem of Petty \cite{MR43:1051} (see Section~\ref{pettysection}) states that $e(X)\leq 2^n$ if the dimension of $X$ is $n$; this then gives $e(\ell_\infty^n)=2^n$.
Kusner \cite{Guy} asked whether $e(\ell_1^n)=2n$ and $e(\ell_p^n)=n+1$ holds for all $1<p<\infty$.
For $p$ an even integer it is easy to show an upper bound of about $pn$ using linear algebra (observation of Galvin; see Smyth \cite{Smyth}).
It is possible to improve this upper bound to about $pn/2$, and in the particular case of $p=4$ to show $e(\ell_4^n)=n+1$ \cite{Swanepoel-AdM} (see Section~\ref{linindsection}).
For $1<p<2$ there are examples known showing $e(\ell_p^n) > n+1$ if $n$ is sufficiently large (depending on $p$) \cite{Swanepoel-AdM} (see Section~\ref{lplowersection}).
On the other hand if $p$ is sufficiently near $2$ (depending on $n$), then $e(\ell_p^n)=n+1$. Smyth \cite{Smyth} gives a quantitative estimate (see Section~\ref{simplexsection}).

The above results are not difficult; to improve the $2^n$ upper bound for $e(\ell_p^n)$ turned out to need deeper results from analysis and linear algebra.
The first breakthrough came when Smyth \cite{Smyth} combined a linear algebra method with the Jackson theorems from approximation theory to prove $e(\ell_p^n) < c_p n^{(p+1)/(p-1)}$ for $1<p<\infty$.
Then Alon and Pudl\'ak \cite{MR1995795} combined Smyth's method with a result on the rank of approximations of the identity matrix (the rank lemma; see Section~\ref{ranksection}) to prove $e(\ell_p^n) < c_p n^{(2p+2)/(2p-1)}$ for $1\leq p<\infty$.
For the $\ell_1$ norm Alon and Pudl\'ak combined the rank lemma with a probabilistic argument (randomized rounding) to prove $e(\ell_1^n)<cn\log n$, and more generally $e(\ell_p^n)< c_p n\log n$ if $p$ is an odd integer.

Finally we mention that it is known that $e(\ell_1^3)=6$ (Bandelt et al. \cite{BCL}) and $e(\ell_1^4)=8$ (Koolen et al.\ \cite{KLS}).
We do not discuss the proofs of these two results.
\begin{problem}[Kusner]
Prove (or disprove) that $e(\ell_1^n)=2n$ for $n\geq 5$.
\end{problem}
\begin{problem}[Kusner]
Prove (or disprove) that $e(\ell_p^n)=n+1$ for $p>2$ ($p\neq 4$).
\end{problem}
\begin{problem}
Prove (or disprove) that $e(\ell_p^n) < cn$ for $1<p<2$.
\end{problem}

\section{General upper bounds}\label{uppersection}
We now need the following notions from convex geometry.
A \emph{polytope} is a convex body that is the convex hull of finitely many points.
A \emph{face} of a convex body $C$ is a set $F\subseteq C$ satisfying the following property:
\begin{quote}
If a segment $\va\vb\subseteq C$ intersects $S$, then $\va\vb\subseteq S$.
\end{quote}
The following are well-known (and easily proved) facts from convex geometry:
\begin{enumerate}
\item Every point of a convex body $C$ is contained in the relative interior of a unique face of $C$.
\item The faces of a polytope forms a ranked finite lattice, with dimension as rank function.
\item The $1$-dimensional faces are called \emph{vertices}, and the $(n-1)$-dimensional faces \emph{facets}.
\item Any $(n-2)$-dimensional face is contained in exactly two facets.
\end{enumerate}
A \emph{parallelotope} is a translate of the polytope $\{\sum_{i=1}^n\lambda_i\va_i:-1\leq\lambda\leq 1\}$ where $\va_1,\dots,\va_n$ are linearly independent.
It is easily seen that if a Minkowski space has a parallelotope as unit ball, then it is isometric to $\ell_\infty^n$.

\subsection{The theorem of Petty and Soltan}\label{pettysection}
The following theorem was proved by Petty \cite{MR43:1051}, and later also by P.~Soltan \cite{MR52:4127}.
A different proof, using the isodiametric inequality, is given in \cite{MR93d:52009}.
\begin{theorem}[Petty \cite{MR43:1051} \& Soltan \cite{MR52:4127}]\label{Petty1}
For any $n$-dimensional normed space $X$, $e(X)\leq 2^n$ with equality iff $X$ is isometric to $\ell_\infty^n$, in which case any equilateral set of size $2^n$ is the vertex set of some ball (which is then a parallelotope).
\end{theorem}
\begin{proof}
Let $\vp_1,\dots,\vp_m\in X$ be $1$-equilateral, where $m=e(X)$.
Let $P=\conv\{\vp_1,\dots,\vp_m\}$ be the convex hull of the points $\vp_i$.
If $P$ is not full-dimensional we may use induction on $n$ to obtain $m\leq 2^{n-1}<2^n$ (and the cases $n=0$ and $n=1$ are trivial).

Thus without loss of generality $P$ is full-dimensional.
For each $i=1,\dots,m$, let $P_i=\frac12(P+\vp_i)$.
Note the following easily proved facts:
\begin{enumerate}
\item Each $P_i\subseteq P$.
\item Each $P_i\subseteq B(\vp_i,\frac12)$.
\item For any $i\neq j$, $P_i\cap P_j$ has empty interior.
\end{enumerate}
Thus the $P_i$ \emph{pack} $P$, and we now calculate volumes:
\[\sum_{i=1}^m\vol{P_i} = \vol{\bigcup_{i=1}^m P_i}\leq\vol{P}.\]
Since $\vol{P_i}=(\frac12)^n\vol{P}$, we obtain $m\leq 2^n$.

If equality holds, $P$ is a polytope in $\R^n$ that is \emph{perfectly packed} or \emph{tiled} by $2^n$ translates of $\frac12 P$.
By the following lemma $P$ must be a parallelotope $\vv+\{\sum_{i=1}^n\lambda_i\va_i:-1\leq\lambda\leq 1\}$ with vertex set $\{\vp_1,\dots,\vp_m\}= \vv+\{\sum_{i=1}^n\epsi_i\va_i: \epsi_i=\pm 1\}$.
Since $\vp_i-\vp_j$ is a unit vector for $i\neq j$, it follows that $\sum_{i=1}^n 2\epsi_i\va_i$ is a unit vector for all choices of $\epsi_i\in\{-1,0,1\}$ with not all $\epsi_i=0$.
It is then easy to see that the unit ball is the parallelotope $\{\sum_{i=1}^n\lambda_i\va_i : -1\leq\lambda\leq1\}$ (exercise), and the linear map sending $\va_i\mapsto\ve_i$ is an isometry.
\end{proof}

\begin{lemma}[Groemer \cite{MR23:A2129}]\label{Groemer}
Let $P$ be the convex hull of $\vv_1,\dots,\vv_{2^n}\in\R^n$.
Suppose that $P=\bigcup_{i=1}^{2^n}\frac12(P+\vv_i)$.
Then $P$ is a parallelotope with vertices $\vv_1,\dots,\vv_{2^n}$.
\end{lemma}
\begin{proof}
Exercise.
\end{proof}

\subsection{Strictly convex norms}
A normed space is strictly convex if $\norm{\vx+\vy}<\norm{\vx}+\norm{\vy}$ for all linearly independent $\vx,\vy$, or equivalently, if the boundary of the uniy ball does not contain a line segment.
If $X$ is a strictly convex Minkowski space, then $e(X)\leq 2^n-1$ by Theorem~\ref{Petty1}, but as far as I know there is no better upper bound for $n\geq 4$.
\begin{problem}
If $X$ is an $n$-dimensional strictly convex Minkowski space, $n\geq 4$, prove that $e(X)\leq 2^n-2$.
\end{problem}
In Section~\ref{slgms} we'll see that $e(X)=3$ when $\dim X=2$, and for $\dim X=3$, Petty \cite{MR43:1051} deduced $e(X)\leq 5$ from a result of Gr\"unbaum \cite{MR28:2480}.
\begin{conjecture}[F\"uredi, Lagarias, Morgan \cite{MR93d:52009}]
There exists some constant $\epsi>0$ such that if $X$ is a strictly convex $n$-dimensional Minkowski space then $e(X)\leq (2-\epsi)^n$.
\end{conjecture}

There exist strictly convex spaces with equilateral sets of size at least exponential in the dimension.
\begin{theorem}[F\"uredi, Lagarias, Morgan \cite{MR93d:52009}]\label{FLM}
For each $n\geq 271$ there exists an $n$-di\-men\-sio\-nal strictly convex $X$ with an equilateral set of size at least $1.058^n$.
\end{theorem}
In \cite{MR93d:52009} the lower bound $1.02^n$ was obtained; however, their argument easily gives the bound as stated here.
They constructed a norm $\norm{\cdot}$ on $\R^n$ with the property that $\norm{\cdot}-\epsi\norm{\cdot}_2$ is still a norm, and asked whether the norm can in addition be $C^\infty$ on $\R^n\setminus\{0\}$.
In our construction the norm will lack the above property, but it will be $C^\infty$.
However, most likely the norm in the above theorem can be made to have the above property and be $C^\infty$.

We need the following very nice high-dimensional phenomenon.
It is a special case of the Johnson-Lindenstrauss flattening lemma \cite{MR86a:46018}, although this special case is essentially the Gilbert-Varshamov lower bound for binary codes (see \cite{vL}.)
\begin{lemma}\label{aorth}
For each $\delta>0$ there exist $\epsi = \epsi(\delta)>0$ and $n_0=n_0(\delta)\geq 1$ such that for all $n\geq n_0$ there exist $m > (1+\epsi)^n$ vectors $\vw_{1},\dots,\vw_{m}\in\R^n$ satisfying 
\begin{equation}\label{iprrel}
\begin{cases}
\ipr{\vw_i}{\vw_i} = 1 & \text{ for all } i, \\
\abs{\ipr{\vw_i}{\vw_j}} < \delta & \text{ for all distinct } i,j.
\end{cases}
\end{equation}
We may take $\epsi=\delta^2/2$ and $n_0\geq (120\log2) /(25\delta^4-\delta^6)$.
\end{lemma}
We need the Chernoff inequality (proved in e.g.\ \cite[Theorem~1.4.5]{vL}).
The {\em binary entropy function} is defined by $H(0) = H(1):= 0$ and
\[
H(x):= -x\binlog x - (1-x)\binlog (1-x), \text{ for }0<x<1.
\]
\begin{lemma}
For any $\epsi>0$ and $n\in\N$,
\begin{equation}\label{Chernoff}
\sum_{0\leq k\leq \epsi n}\binom{n}{k}\leq 2^{nH(\epsi)}.
\end{equation}
\end{lemma}

\begin{proof}[Proof of Lemma~\ref{aorth}]
For $n\geq 1$, define a graph on $\{-1,1\}^n$ as vertex set by joining two points if the number of coordinates in which they differ is $\leq \frac{1}{2}(1-\delta)n$ or $\geq\frac{1}{2}(1+\delta)n$.
This graph is regular of degree 
$$m:=\sum_{0\leq k\leq\frac{1}{2}(1-\delta)n}\binom{n}{k}+\sum_{\frac{1}{2}(1+\delta)n\leq k\leq n}\binom{n}{k}\leq -1+2^{nH((1-\delta)/2)+1},$$
by \eqref{Chernoff}.
If $\vx$ and $\vy$ are not connected in this graph, then $n^{-1/2}\vx$ and $n^{-1/2}\vy$ are unit vectors satisfying $\abs{\ipr{\vx}{\vy}}<\delta$.
We therefore have to show that the graph contains an \emph{independent set} (a set without an edge between any two vertices) of at least $(1+\epsilon)^n$ points.
We arbitrarily choose a vertex, delete it and its $m$ neighbours, choose a remaining vertex, delete it and its at most $m$ neighbours, and continue until nothing remains.
Thus we obtain an independent set of size at least
\[2^n/(m+1) \leq 2^{n(1-H((1-\delta)/2))-1}\]
vertices, which is greater than $(1+\delta^2/2)^n$ if
\[n f(\delta)>1,\]
where $f(\delta)=1-H((1-\delta)/2)-\binlog(1+\delta^2/2)$.
But this is true for sufficiently large $n$, since it is easily seen that $f(\delta)>0$ for $0<\delta<1$.
That $n\geq (120\log2) /(25\delta^4-\delta^6)$ is sufficient can be seen by calculating the Taylor expansion of $f(\delta)$.
\end{proof}
\subsubsection*{Remarks}
\begin{enumerate}
\item The proof shows that if $M(n,\delta)$ is the maximum number of unit vectors in $\R^n$ satisfying \eqref{iprrel}, then
\[ \lim_{n\to\infty}\frac{1}{n}\binlog M(n,\delta) \geq 1-H\Big(\frac{1-\delta}{2}\Big) = \frac{1}{\log 2}\Big(\frac{\delta^2}{2}+\frac{\delta^4}{12}+\dots\Big).\]
If one repeatedly choose unit vectors on the Euclidean unit sphere $\Sphere^{n-1}$ and delete spherical caps around them, then, as shown by Wyner \cite[pp.~1089--1092]{Wyner} (slightly improving an earlier estimate of Shannon \cite{Shannon}),
\[ \lim_{n\to\infty}\frac{1}{n}\binlog M(n,\delta) \geq \binlog (1-\delta^2)^{-1/2} = \frac{1}{\log 2}\Big(\frac{\delta^2}{2}+\frac{\delta^4}{4}+\dots\Big),\]
which is the same up to first order, but has a better higher order term.
Thus the above greedy argument is not really improved by choosing vectors from the Euclidean unit sphere.
\item Recently, the Gilbert-Varshamov lower bound for binary codes has been slightly improved \cite{JV}, which will improve the lower bound in Lemma~\ref{aorth} to $c(\delta)n(1+\delta^2/2)^n$.
\item If we prove this lemma by choosing vectors in $\{\pm 1\}$ randomly, and estimating the probabilities in the obvious way, then we obtain a slightly worse dependence of $\epsi=\delta^2/4$.
\end{enumerate}

\begin{proof}[Proof of Theorem~\ref{FLM}]
By Lemma~\ref{aorth} there exists $\vw_1,\dots,\vw_m$, $m\geq 1.058^n$, satisfying \eqref{iprrel} with $\delta=1/3$.

We now show that $P=\conv\{\vw_i-\vw_j: 1\leq i,j\leq m\}$ is a centrally symmetric polytope with $m(m-1)$ vertices $\vw_i-\vw_j$, $i\neq j$.
It is sufficient to note that the hyperplane \[\{\vx\in\R^n:\ipr{\vx}{\vw_i-\vw_j} = \ipr{\vw_i-\vw_j}{\vw_i-\vw_j}\}\] is a supporting hyperplane of $P$ passing only through $\vw_i-\vw_j$.
This follows from
\begin{eqnarray*}
&& \ipr{\vw_j-\vw_i}{\vw_i-\vw_j}, \pm\ipr{\vw_i-\vw_k}{\vw_i-\vw_j},\\
&& \pm\ipr{\vw_j-\vw_k}{\vw_i-\vw_j}, \pm\ipr{\vw_k-\vw_\ell}{\vw_i-\vw_j}\\
&\leq& \ipr{\vw_i-\vw_j}{\vw_i-\vw_j} \text{ for distinct $i,j,k,\ell$,}
\end{eqnarray*}
which follow easily from \eqref{iprrel} and $\delta=1/3$.

With $P$ as a unit ball we thus already have $\{\vw_i:i=1,\dots,m\}$ as an equilateral set.
We now need to find a strictly convex centrally symmetric body with the vertices of $P$ on its boundary.
This is provided by the next lemma.
\end{proof}

The following lemma is obvious, but we provide a simple proof.
\begin{lemma}\label{smoothexist}
Let $S$ be the vertex set of a centrally symmetric polytope in $\R^n$.
Then there is a smooth, strictly convex norm $\norm{\cdot}$ on $\R^n$ such that $\norm{\vx}=1$ for all $\vx\in S$.
\end{lemma}
\begin{proof}
Let $S=\{\pm \vx_1, \dots, \pm \vx_m\}$.
For each $\vx_i$, choose $\vy_i\in\R^n$ such that $\ipr{\vy_i}{\vx_i}=1$ and $\abs{\ipr{\vy_i}{\vx_j}}<1$ for all $i\neq j$.
The required norm will be
$$\norm{\vx} = \Big(\sum_{j=1}^m\lambda_j\abs{\ipr{\vy_j}{\vx}}^p\Big)^{1/p}$$
for suitably chosen $\lambda_j>0$ and $1<p<\infty$, i.e., we want to imbed $S$ into the unit sphere of $\ell_p^m$ for some $p$.

Choose $p$ large enough such that
$$\abs{\ipr{\vy_j}{\vx_i}}^p<\frac{1}{2m} \text{ for all } j\neq i.$$
Consider the matrix $A=\begin{bmatrix}a_{ij}\end{bmatrix}_{i,j=1}^m$ with
$$a_{ij}=\begin{cases}\abs{\ipr{\vy_j}{\vx_i}}^p & \text{ if } i\neq j, \\ 0 & \text{ if } i=j.\end{cases}$$
Considering $A$ as a linear transformation on $\ell_\infty^m$, we have $\norm{A}<\half$.
We want to find a vector $\vv=(\lambda_1,\dots,\lambda_m)^\Tr > 0$ such that
$$(I+A)\vv = (1,\dots,1)^\Tr,$$
since then we would have $\norm{\vx_i}=1$ for all $i$.
However, since $\norm{A}<\half$, $I+A$ is invertible, and then necessarily $\vv=(I+A)^{-1}(1,\dots,1)^\Tr$.
Also,
\begin{align*}
\norm{\vv-(1,\dots,1)^\Tr}_\infty & = \norm{(I+A)^{-1}(1,\dots,1)^\Tr-(1,\dots,1)^\Tr}_\infty\\
& \leq \norm{(I+A)^{-1}-I}\norm{(1,\dots,1)}_\infty\\
& = \norm{\sum_{i=1}^\infty A^i} \leq \sum_{i=1}^\infty\norm{A}^i < 1.
\end{align*}
Thus, $\lambda_i>0$ for all $i$.
\end{proof}

Note that in the above lemma we may choose $p$ to be an even integer which gives a norm that is $C^\infty$ on $\R^n\setminus\{\vo\}$.

\section{Simple lower bounds}
We now consider the problem of finding equilateral sets in a general Minkowski space, thus providing a lower bound to $e(X)$.

\subsection{General Minkowski spaces}\label{slgms}
\begin{proposition}\label{twod}
If $\dim X\geq 2$ then $e(X)\geq 3$.
\end{proposition}
\begin{proof}
Exercise. (Hint: Euclid Book I Proposition 1.)
\end{proof}

The above proposition combined with Theorem~\ref{Petty1} gives the following
\begin{corollary}\label{e2d}
If $\dim X=2$ then
\begin{equation*}
e(X)=\begin{cases} 4 & \text{if the unit ball of $X$ is a parallelogram},\\
                   3 & \text{if the unit ball of $X$ is not a parallelogram.}
\end{cases}
\end{equation*}
\end{corollary}

In general one would hope for the following conjecture, stated in \cite{MR25:1492, MR93h:53012, MR43:1051, MR97f:52001}.
\begin{conjecture}\label{equilconj}
If $\dim X=n$ then $e(X)\geq n+1$.
\end{conjecture}
As seen above this is simple for $n=2$, and the next theorem shows its truth for $n=3$.
However, for each $n\geq 4$ this is open, and the best that is known is the theorem of Brass and Dekster giving $e(X)\geq c(\log n)^{1/3}$, discussed in Section~\ref{Brasssection}.

\begin{theorem}[Petty \cite{MR43:1051}]\label{Petty2}
If $\dim X\geq 3$ then any equilateral set of size $3$ may be extended to an equilateral set of size $4$.
\end{theorem}

The proof needs the following technical but simple lemma.

\begin{lemma}\label{monotone}
Let $S_2$ be the unit circle (i.e.\ boundary of the unit ball) of a $2$-dimensional Minkowski space.
Fix any $\vu\in S_2$.
Let $\vf:[0,1]\to S_2$ be a parametrization of an arc of $S_2$ from $\vu$ to $-\vu$.
Fix $\vp=\lambda\vu$, $\lambda>0$.
Then $t\mapsto\norm{\vp-\vf(t)}$ is strictly increasing before it reaches the value $\norm{\vp+\vu}$, i.e.\ on the interval $[0,t_0]$, where $t_0=\min\{t:\norm{\vp-\vf(t)}<\norm{\vp+\vu}\}$.
\end{lemma}
\begin{proof}
Exercise.
\end{proof}

\begin{proof}[Proof of Theorem~\ref{Petty2}]
Let $\va,\vb,\vc$ be $1$-equilateral in $X$.
Without loss of generality $\va=\vo$.
Let $X_2=\Span\{\vb,\vc\}$, and let $X_3$ be any $3$-dimensional subspace of $X$ containing $X_2$.
Let $S_2$ be the unit circle of $X_2$, and let $S_3$ be the unit sphere of $X_3$.

Consider the mapping $\vf:S_3\to\R^2$ defined by $\vf(x)=(\norm{\vx-\vb},\norm{\vx-\vc})$.
We have to show that $(1,1)$ is in the range of $\vf$.
If we can show that the restriction of $\vf$ to $S_2$ encircles $(1,1)$, then the theorem follows from the following topological argument.

$S_2$ can be contracted to a point on $S_3$, since $S_3$ is simply connected; thus $\vf(S_2)$ can also be contracted to a point; if $(1,1)$ is in the interior of the Jordan curve $\vf(S_2)$, then at some stage of the contraction the curve must pass through $(1,1)$, since $\R^2\setminus\{(1,1)\}$ is not simply connected.

Thus we assume without loss of generality that $(1,1)\notin \vf(S_2)$, and we now follow the curve $\vf(S_2)$.
We start at $\vf(\vb)=(0,1)$ and go to $\vf(\vc)=(1,0)$.
By Lemma~\ref{monotone} the $x$-coordinate increases strictly while the $y$-coordinate decreases strictly.
Therefore the arc of the curve strictly between $(0,1)$ and $(1,0)$ is contained in the square $0<x,y<1$.
Then from $(1,0)$ to $\vf(\vc-\vb)=(\norm{\vc-2\vb},1)$ the $y$-coordinate increases strictly, and the $x$-coordinate also increases strictly at least in a neighbourhood of $(1,0)$ if $\norm{\vc-2\vb}>1$.
However, by the triangle inequality we have $\norm{\vc-2\vb}\geq1$, and if $\norm{\vc-2\vb}=1$ then we already have found $(1,1)$ to be in the range of $\vf$, which we assumed not to happen.

From $\vf(\vc-\vb)$ to $\vf(-\vb)=(2,\norm{\vb+\vc})$ the $x$-coordinate increases, and from $\vf(-\vb)$ to $\vf(-\vc)=(\norm{\vb+\vc},2)$ the $y$-coordinate increases.
From $\vf(-\vc)$ to $\vf(\vb-\vc)=(1,\norm{\vb-2\vc})$ the $x$-coordinate decreases strictly, and the $y$-coordinate decreases.
Again by the triangle inequality and assumption we have $\norm{\vb-2\vc} > 1$.
Finally, from $\vf(\vb-2\vc)$ to $\vf(\vb)$ the $x$-coordinate decreases strictly, and the $y$-coordinate decreases.

It follows that $\vf(S_2)$ encircles $(1,1)$, and the theorem is proved.
\end{proof}

Note that if the fourth point which extends the equilateral set is in the plane $X_2$ (i.e.\ if $(1,1)\in \vf(S_2)$ in the above proof), then $X_2$ is isometric to $\ell_\infty^2$, i.e.\ $S_2$ is a parallelogram (by Corollary~\ref{e2d}).
However, it is impossible for all $2$-dimensional sections of the unit ball through the origin to be parallelograms.
Thus we may always choose an $X_2$ so that the fourth point is outside $X_2$ to obtain a non-coplanar equilateral set of four points.

One would hope that ideas similar to the above proof would help to prove Conjecture~\ref{equilconj}.
However, the above proof cannot be na\"\i vely generalized, as the following example of Petty \cite{MR43:1051} shows.

Define the following norm on $\R^n$:
\[ \norm{(x_1,\dots,x_n)} := \abs{x_1} + \sqrt{x_2^2+\dots+x_n^2}.\]
Thus $\norm{\cdot}$ is the $\ell_1$-sum of $\R$ and $\ell_2^{n-1}$: $(\R^n,\norm{\cdot})=\R\oplus_1\ell_2^{n-1}$.
Its unit ball is a double cone over an $(n-1)$-dimensional Euclidean ball:
\[ B_{\norm{\cdot}} = \conv(B_2^{n-1}\cup\{\pm\ve_1\}).\]

For any $n\geq 2$, $(\R^n,\norm{\cdot})$ contains a maximal equilateral set of four points.
In fact the largest size of an equilateral set containing $\pm\ve_1$ is $4$.
First of all note that
\[ \norm{\vx-\ve_1}=\norm{\vx+\ve_1} \iff x_1=0.\]
Secondly, since $\norm{\ve_1-(-\ve_1)}=2$, we want $\norm{\vx\pm\ve_1}=2$.
Since then $x_1=0$, we obtain that $1+\sqrt{x_2^2+\dots+x_n^2}=2$ must hold, which gives that $\vx\in\Sphere^{n-2}$.
Clearly there are at most two points on the Euclidean $(n-2)$-dimensional unit sphere at distance $2$ (and two such points must be antipodal).

On the other hand any $n$ equilateral points in the $(n-1)$-dimensional Euclidean subspace $x_1=0$ can always be extended to $n+1$ equilateral points.
For example, if we choose the $n$ equilateral points on $\Sphere^{n-2}$ then $\lambda\ve_1$ will be the $(n+1)$st point, for some appropriate $\lambda\in\R$.

\begin{problem}
Calculate $e(X)$ for the Minkowski space in the above example.
\end{problem}

\subsection{$\ell_p^n$ with $1<p<2$}\label{lplowersection}
We sketch the proof that $e(\ell_p^n)>n+1$ for any $1<p<2$ with $n$ sufficiently large (depending on $p$).

Note that the four $\pm1$-vectors in $\R^3$ with an even number of $-1$'s, i.e.,
\begin{equation}\label{tetra}
 \begin{bmatrix}  1\\ 1\\ 1\end{bmatrix}, 
 \begin{bmatrix}  1\\-1\\-1\end{bmatrix}, 
 \begin{bmatrix} -1\\ 1\\-1\end{bmatrix}, 
 \begin{bmatrix} -1\\-1\\ 1\end{bmatrix},
\end{equation}
is an equilateral set for any $p$-norm.
Thus if we consider $\R^6$ to be the direct sum $\R^3\oplus\R^3$ and put the above four vectors in each copy of $\R^3$, we obtain the following $8$ vectors in $\R^6$:
\begin{equation}
 \begin{bmatrix}  1\\ 1\\ 1\\ 0\\ 0\\ 0\end{bmatrix}, 
 \begin{bmatrix}  1\\-1\\-1\\ 0\\ 0\\ 0\end{bmatrix}, 
 \begin{bmatrix} -1\\ 1\\-1\\ 0\\ 0\\ 0\end{bmatrix}, 
 \begin{bmatrix} -1\\-1\\ 1\\ 0\\ 0\\ 0\end{bmatrix},
 \begin{bmatrix}  0\\ 0\\ 0\\ 1\\ 1\\ 1\end{bmatrix}, 
 \begin{bmatrix}  0\\ 0\\ 0\\ 1\\-1\\-1\end{bmatrix}, 
 \begin{bmatrix}  0\\ 0\\ 0\\-1\\ 1\\-1\end{bmatrix}, 
 \begin{bmatrix}  0\\ 0\\ 0\\-1\\-1\\ 1\end{bmatrix}.
\end{equation}
Then it is easily seen that these $8$ vectors are equilateral iff $2\cdot2^p = 6\iff p=\log3/\log2$.

The combinatorial reason why this worked is that any two of the vectors in \eqref{tetra} differ in exactly two positions.
For a generalization we use Hadamard matrices.
An \emph{Hadamard matrix} is an $n\times n$ matrix $H$ with all entries $\pm 1$ and with orthogonal columns.
By multiplying some of the columns of an Hadamard matrix one obtains an Hadamard matrix with the first row containing only $1$'s.
If we remove this first row, we obtain $n$ column vectors in $\R^{n-1}$ with $\pm1$ entries such that any two differ in exactly $n/2$ positions.
The vectors in \eqref{tetra} was obtained in this way from a $4\times 4$ Hadamard matrix.
We may now construct $2n$ vectors in $\R^{2n-2}=\R^{n-1}\oplus\R^{n-1}$ as before that will be equilateral iff
\[\frac{n}{2}2^p=2n-2 \iff p=2+\binlog\frac{n-1}{n}.\]

The only problem remaining is to find Hadamard matrices.
It is easy to see that for an $n\times n$ Hadamard matrix to exist, $n$ must be divisible by $4$.
It is an open problem whether the converse is true.
However, the following well-known construction gives a Hadamard matrix of order $2^k$:
\[ H_0 = \begin{bmatrix} 1 \end{bmatrix}, \quad H_{k+1}=\begin{bmatrix} H_k & H_k\\ H_k & -H_k \end{bmatrix},\; k\geq 0.\]
Thus we have $e(\ell_p^n)>n+1$ for $p$ arbitrarily close to $2$.

With a more involved construction the following can be shown:
\begin{theorem}
For any $1<p<2$ and $n\geq 1$, let
\[k=\left\lceil\frac{\log(1-2^{p-2})^{-1}}{\log 2}\right\rceil-1.\]
Then
\[e(\ell_p^n)\geq \left\lfloor\frac{2^{k+1}}{2^{k+1}-1}n\right\rfloor.\]
In particular, if $n\geq 2^{k+2}-2$ then $e(\ell_p^n)>n+1$.
\end{theorem}
\begin{corollary}\label{Swancor}
$e(\ell_p^n)>n+1$ if $p< 2- \frac{1+o(1)}{(\log 2) n}$.
\end{corollary}
Compare this with Corollary~\ref{Smythcor} in Section~\ref{simplexsection}.
The smallest dimension for which the theorem gives an example of $e(\ell_p^n)>n+1$ is $n=6$ (and $1<p\leq\log3/\log 2$).
In \cite{Swanepoel-AdM} there is also a $4$-dimensional example (with  $1<p\leq\log\frac52/\log 2$).

In the next two sections we show that $e(X)$ goes to infinity with the dimension of $X$ (Brass-Dekster), although the best known lower bound for $e(X)$ is much smaller than linear in $n$.
There are three ingredients to the proof: The Cayley-Menger theory of the embedding of metric spaces into Euclidean space, Dvoretzky's theorem, and the Brouwer fixed point theorem.

\section{Cayley-Menger Theory}
In this section we discuss a fragment of the Cayley-Menger theory by giving necessary and sufficient conditions for a metric space of size $n+1$ to be embeddable as an affinely independent set in $\ell_2^n$.
The general theory is by Menger \cite{Menger}; see also Blumenthal \cite[\S40]{MR42:3678}.

\subsection{The Cayley-Menger determinant}
Let a metric space on $n+1$ points $\vp_0,\dots,\vp_n$ be given, with distances $\rho_{ij}=d(\vp_i,\vp_j)$.
We now derive a necessary condition for $\vp_0,\dots,\vp_n$ to be isometric to an affinely independent subset of $\ell_2^n$, in terms of $\rho_{ij}$.

So we assume that $\{\vp_0,\dots,\vp_n\}$ is an affinely independent set in $\ell_2^n$ with $\norm{\vp_i-\vp_j}_2=\rho_{ij}$.
We write the coordinates of $\vp_i$ as $(p_i^{(1)},p_i^{(2)},\dots,p_i^{(n)})$.
Let $\Delta=\conv\{\vp_0,\dots,\vp_n\}$ denote the simplex with vertices $\vp_i$.
Then
\begin{eqnarray*}
\pm\vol{\Delta} & = & \frac{1}{n!}
\begin{vmatrix}
p_0^{(1)} & p_1^{(1)} & \cdots & p_n^{(1)}\\
p_0^{(2)} & p_1^{(2)} & \cdots & p_n^{(2)}\\
\vdots & \vdots & \cdots & \vdots\\
p_0^{(n)} & p_1^{(n)} & \cdots & p_n^{(n)}\\
1 & 1 & \cdots & 1
\end{vmatrix}\\
&=& \frac{1}{n!}
\begin{vmatrix}
\vp_0 & \vp_1 & \cdots & \vp_n\\
1 & 1 & \cdots & 1
\end{vmatrix}.
\end{eqnarray*}
Squaring (and adding an extra row and column),
\begin{eqnarray*}
\vol{\Delta}^2 &=& \frac{1}{(n!)^2}\det
\begin{bmatrix}
\vp_0^\Tr & 1 & 0\\
\vp_1^\Tr & 1 & 0\\
\vdots & \vdots & \vdots\\
\vp_n^\Tr & 1 & 0\\
\vo^\Tr & 0 & 1
\end{bmatrix}
\begin{bmatrix}
\vp_0 & \cdots & \vp_n & \vo \\
1 & \cdots & 1 & 0 \\
0 & \cdots & 0 & 1
\end{bmatrix}\\
&=& -\frac{1}{(n!)^2}\det
\begin{bmatrix}
\vp_0^\Tr & 0 & 1\\
\vp_1^\Tr & 0 & 1\\
\vdots & \vdots & \vdots\\
\vp_n^\Tr & 0 & 1\\
\vo^\Tr & 1 & 0
\end{bmatrix}
\begin{bmatrix}
\vp_0 & \cdots & \vp_n & \vo \\
1 & \cdots & 1 & 0 \\
0 & \cdots & 0 & 1
\end{bmatrix} \quad
\parbox{5cm}{(interchange last two \\ \mbox{}\quad columns of first matrix)}\\
&=& -\frac{1}{(n!)^2}
\begin{vmatrix}
\ipr{\vp_0}{\vp_0} & \ipr{\vp_0}{\vp_1} & \cdots & \ipr{\vp_0}{\vp_n} & 1\\
\ipr{\vp_1}{\vp_0} & \ipr{\vp_1}{\vp_1} & \cdots & \ipr{\vp_1}{\vp_n} & 1\\
\vdots & \vdots & \cdots & \vdots\\
\ipr{\vp_n}{\vp_0} & \ipr{\vp_n}{\vp_1} & \cdots & \ipr{\vp_n}{\vp_n} & 1\\
1 & 1 & \cdots & 1 & 0
\end{vmatrix}\\
&=& -\frac{1}{2^{n+1}(n!)^2}
\begin{vmatrix}
0 & -\rho_{01}^2 & \cdots & -\rho_{0n}^2 & 1 \\
-\rho_{10}^2 & 0 & \ddots & \vdots & \vdots \\
\vdots & \ddots & \ddots & \vdots & \vdots\\
-\rho_{n0}^2 & \cdots & -\rho_{n,n-1}^2 & 0 & 1 \\
1 & \cdots & \cdots & \cdots & 0
\end{vmatrix} \qquad\parbox{4cm}{(use $2\ipr{\vp_i}{\vp_j}=$\\$\ipr{\vp_i}{\vp_i}+\ipr{\vp_j}{\vp_j}-\rho_{ij}^2$ and subtract multiples of the last row and column)}\\
&=& \frac{(-1)^{n+1}}{2^{n+1}(n!)^2}\det
\begin{bmatrix}
 & & & 1 \\
 &P& & \vdots \\
 & & & 1 \\
1&\cdots&1&0
\end{bmatrix},
\end{eqnarray*}
where $P=\begin{bmatrix}\rho_{ij}^2\end{bmatrix}$.
We call
\[\CMdet(\vp_0,\dots,\vp_n) :=
\begin{vmatrix}
 & & & 1 \\
 &P& & \vdots \\
 & & & 1 \\
1&\cdots&1&0
\end{vmatrix}\]
the \emph{Cayley-Menger determinant} of the $(n+1)$-point metric space.
We have shown the following.
\begin{proposition}\label{CMnecessity}
If a metric space on $n+1$ points can be embedded into $\ell_2^n$ as an affinely independent set, then its Cayley-Menger determinant has sign $(-1)^{n+1}$.
(Also, if the metric space can be embedded as an affinely dependent set, then its Cayley-Menger determinant is $0$.)
\end{proposition}

Conversely we have
\begin{theorem}[Menger \cite{Menger}]\label{CMsufficiency}
A metric space on $n+1$ points $\vp_0,\dots,\vp_n$ can be embedded into $\ell_2^n$ as an affinely independent set if for each $k=1,\dots,n$, the Cayley-Menger determinant of $\vp_0,\dots,\vp_k$ has sign $(-1)^{k+1}$.
\end{theorem}
\begin{proof}
For $n=1$,
\[ \CMdet(\vp_0,\vp_1) = 
\begin{vmatrix}
0 & \rho_{01}^2 & 1\\
\rho_{10}^2 & 0 & 1\\
1 & 1 & 0
\end{vmatrix}
= \rho_{01}^2+\rho_{10}^2 > 0.\]
Thus the theorem is trivial for $n=1$ (and even more so for $n=0$).

Thus we may assume as inductive hypothesis that the theorem is true for $n-1$.
Thus $\vp_0,\dots,\vp_{n-1}$ can be embedded into $\ell_2^{n-1}\subset\ell_2^n$ as an affinely independent set, say $\vp_0\mapsto\vo$ and $\vp_i\mapsto\vx_i$, $i=1,\dots,n-1$.
Then $\Span\{\vx_1,\dots,\vx_{n-1}\}=\ell_2^{n-1}$ and $\vx_1,\dots,\vx_{n-1}$ are linearly independent.

We have $\norm{\vx_i}_2=\rho_{i0}$, and $\norm{\vx_i-\vx_j}_2=\rho_{ij}$ for $1\leq i,j\leq n-1$, and we have to find an $\vx\notin \ell_2^{n-1}$ satisfying
\begin{equation}\label{req1}\norm{\vx}_2=\rho_{0n}
\end{equation}
and
\begin{equation}\label{req2}
\norm{\vx-\vx_i}_2=\rho_{in}\text{ for all }i=1,\dots,n-1.
\end{equation}
We write $\vx=\vv+\lambda\ve_n$ uniquely, where $\vv\in \ell_2^{n-1}$ and $\lambda\in\R$.
If we square \eqref{req2} and simplify, we find that $\vv$ has to satisfy the following $n-1$ linear equations in $\ell_2^{n-1}$:
\[ 2\ipr{\vx_i}{\vv} = \norm{\vx}_2^2 + \rho_{0i}^2-\rho_{in}^2,\; i=1,\dots,n-1.\]
Since the $\vx_i$ are linearly independent, there is a unique solution $\vv\in\ell_2^{n-1}$ to the following modified equation:
\begin{equation}\label{foundv}
2\ipr{\vx_i}{\vv} = \rho_{0n}^2 + \rho_{0i}^2-\rho_{in}^2,\; i=1,\dots,n-1,
\end{equation}
and we will have found the required $\vx$ once we can satisfy \eqref{req1}, which is equivalent to
\[ \lambda^2=\rho_{0n}^2-\norm{\vv}_2^2.\]
Thus it remains to prove that $\rho_{0n}^2-\norm{\vv}_2^2>0$.
To this end we calculate the Cayley-Menger determinant of $\vo,\vx_1,\dots,\vx_{n-1},\vv$, which is $0$ by Proposition~\ref{CMnecessity}:
\[
\begin{vmatrix}
0 & \rho_{01}^2 & \cdots & \rho_{0,n-1}^2 & \norm{\vv}_2^2 & 1\\
\rho_{10}^2 & 0 & \cdots & \rho_{1,n-1}^2 & \norm{\vv-\vx_1}_2^2 & 1\\
\vdots & \ddots & \ddots & \vdots & \vdots & \vdots\\
\rho_{n-1,0}^2 & \cdots & \cdots & 0 & \norm{\vv-\vx_{n-1}}_2^2 & 1\\
\norm{\vv}_2^2 & \norm{\vv-\vx_1}_2^2 & \cdots & \norm{\vv-\vx_{n-1}}_2^2 & 0 & 1\\
1 & 1 & \cdots & 1 & 1 & 0
\end{vmatrix}
=0.\]
Multiplying out
\begin{eqnarray*}
\norm{\vv-\vx_i}_2^2 & = & \norm{\vv}_2^2-2\ipr{\vx_i}{\vv}+\norm{\vx_i}_2^2\\
& = & \norm{\vv}_2^2 -\rho_{0n}^2+\rho_{in}^2\quad\text{(by \eqref{foundv}),}
\end{eqnarray*}
and using the last row and column to eliminate $\norm{\vv}_2^2-\rho_{0n}^2$, we obtain
\[
\begin{vmatrix}
0 & \rho_{01}^2 & \cdots & \rho_{0,n-1}^2 & \rho_{0n}^2 & 1\\
\rho_{10}^2 & 0 & \cdots & \rho_{1,n-1}^2 & \rho_{1n}^2 & 1\\
\vdots & \ddots & \ddots & \vdots & \vdots & \vdots\\
\rho_{n-1,0}^2 & \cdots & \cdots & 0 & \rho_{n-1,n}^2 & 1\\
\rho_{n,0}^2 & \rho_{n,1}^2 & \cdots & \rho_{n,n-1}^2 & -2(\norm{\vv}_2^2-\rho_{0n}^2) & 1\\
1 & 1 & \cdots & 1 & 1 & 0
\end{vmatrix}
=0.\]
This determinant differs from $\CMdet(\vp_0,\dots,\vp_n)$ only in the second last column.
Thus we may subtract this determinant from $\CMdet(\vp_0,\dots,\vp_n)$, which we are given has sign $(-1)^{n+1}$, to obtain the following determinant of sign $(-1)^{n+1}$:
\begin{eqnarray*}
&&\begin{vmatrix}
0 & \rho_{01}^2 & \cdots & \rho_{0,n-1}^2 & 0 & 1\\
\rho_{10}^2 & 0 & \cdots & \rho_{1,n-1}^2 & 0 & 1\\
\vdots & \ddots & \ddots & \vdots & \vdots & \vdots\\
\rho_{n-1,0}^2 & \cdots & \cdots & 0 & 0 & 1\\
\rho_{n,0}^2 & \rho_{n,1}^2 & \cdots & \rho_{n,n-1}^2 & 2(\norm{\vv}_2^2-\rho_{0n}^2) & 1\\
1 & 1 & \cdots & 1 & 1 & 0
\end{vmatrix}\\
&=& 2(\norm{\vv}_2^2-\rho_{0n}^2)\CMdet(\vp_0,\dots,\vp_{n-1})\quad\text{(expanding along 2nd last column)}.
\end{eqnarray*}
Since it is also given that $\CMdet(\vp_0,\dots,\vp_{n-1})$ has sign $(-1)^{n}$, it follows that $\norm{\vv}_2^2-\rho_{0n}^2<0$.
This gives us a $\lambda\in\R$ such that $x=\vv+\lambda\ve_n$ satisfies \eqref{req1} as well.
The theorem is proved.
\end{proof}

\subsection{Embedding an almost equilateral $n$-simplex into $\ell_2^n$}\label{simplexsection}
We now want to use Theorem~\ref{CMsufficiency} to show that an $n$-simplex with all side lengths close to $1$ can be embedded isometrically into $\ell_2^n$.
Thus we have to find the sign of the determinant of a matrix with $0$ on the diagonal and all off-diagonal entries close to $1$, since this is how the Cayley-Menger determinant will look like.

We start with the following simple observations.
Let $J$ be the $n\times n$ matrix with all entries equal to $1$.
\begin{lemma}\label{detlemma}
Let $n\geq 1$.
\begin{enumerate}
\item The characteristic polynomial of $J$ is $(-1)^n(x-n)x^{n-1}$.
\item Therefore the characteristic polynomial of $J-I$ is $(-1)^n(x-n+1)(x+1)^{n-1}$, and $\det(J-I)=(-1)^{n-1}(n-1)$.
\item Therefore $J-I$ is invertible, and its inverse has eigenvalues $-1$ with multiplicity $n-1$, and $1/(n-1)$.
\end{enumerate}
\end{lemma}

\begin{proposition}\label{detsignlemma}
Let $A=\begin{bmatrix}a_{ij}\end{bmatrix}$ be an $(n+2)\times (n+2)$ matrix with $a_{ii}=0$, $\abs{1-a_{ij}}\leq 1/(n-\frac32)$ for $1\leq i\neq j\leq n+1$, and $a_{i,n+2}=a_{n+2,i}=1$ for $i=1,\dots,n$.
Then $\det(A)$ has sign $(-1)^{n+1}$.
\end{proposition}
\begin{proof}
We first prove that $A$ is invertible.
Let $E=\begin{bmatrix}\epsi_{ij}\end{bmatrix}:=J-I-A$.
Then $\epsi_{ii}=0$ and $\epsi_{i,n+2}=\epsi_{n+2,i}=0$ for all $i$, and $\abs{\epsi_{ij}}\leq1/(n-\frac32)$ for all $1\leq i\neq j\leq n+1$.
Since $A=J-I-E=(J-I)(I-(J-I)^{-1}E)$, it is sufficient to prove that $\norm{(J-I)^{-1}E}_2<1$, where $\norm{\cdot}_2$ is the operator norm on $L(\ell_2^n\to\ell_2^n)$.
Firstly, since the operator norm is the maximum modulus of the eigenvalues, we have by Lemma~\ref{detlemma} that $\norm{(J-I)^{-1}}_2=1$.
 Secondly, by the Cauchy-Schwartz inequality, $\norm{E}_2\leq\sqrt{\sum_{i,j=1}^{n+2}\epsi_{ij}^2}<1$.
Thus $A$ is invertible.

To show that $\det A$ has sign $(-1)^{n+1}$, observe that we may join $J-I$ to $A$ with a curve $A_t$, $0\leq t\leq 1$, where each $A_t$ satisfies the hypothesis of the proposition.
Thus we know from what we have just proved that each $A_t$ is invertible.
By continuity considerations the sign of $\det A_t$ does not change.
Therefore, $A$ and $J-I$ have the same sign, which is $(-1)^{n+1}$ by Lemma~\ref{detlemma}.
\end{proof}

\begin{theorem}[Dekster \& Wilker \cite{MR88k:51035}]\label{simplexemb}
Let $\vp_0,\dots,\vp_n$ be a metric space with distances satisfying
\[1\leq d(\vp_i,\vp_j)\leq 1+\frac{1}{n+1}.\]
Then $\vp_0,\dots,\vp_n$ can be embedded into $\ell_2^n$, and any such embedding must be an affinely independent set.
\end{theorem}
\begin{proof}
We first scale the distances by some $\alpha>0$ in such a way that the given inequality is transformed into
\[ 1-\epsi \leq (\alpha d(\vp_i,\vp_j))^2\leq 1+\epsi\]
for some $\epsi>0$.
Thus $\epsi$ and $\alpha$ must satisfy
\[ \frac{1-\epsi}{\alpha^2}=1\text{ and } \frac{1+\epsi}{\alpha^2}=\left(\frac{n+2}{n+1}\right)^2.\]
Eliminating $\alpha$ we find
\begin{equation}\label{epsilon}
\frac{1+\epsi}{1-\epsi}=\left(\frac{n+2}{n+1}\right)^2.
\end{equation}
We want to apply Theorem~\ref{CMsufficiency}.
Since the CM determinant of $\vp_0,\dots,\vp_k$ is a $(k+2)\times(k+2)$ matrix with $k\leq n$, to show that it has sign $(-1)^{k+1}$, it is sufficient to prove that
\[\epsi\leq\frac{1}{n+\frac32},\]
according to Proposition~\ref{detsignlemma}.
By \eqref{epsilon} we have
\begin{eqnarray*}
\epsi&=& \frac{\left(\frac{n+2}{n+1}\right)^2-1}{\left(\frac{n+2}{n+1}\right)^2+1}\\
&=& \frac{(n+2)^2-(n+1)^2}{(n+2)^2+(n+1)^2} = \frac{2n+3}{2n^2+6n+5}\\
&<& \frac{2n+3}{2n^2+6n+\frac92} = \frac{2}{2n+3}.
\end{eqnarray*}
By Proposition~\ref{CMnecessity} no embedding can be affinely dependent.
\end{proof}

Dekster and Wilker found the smallest $\gamma_n>1$ such that whenever
\[1\leq d(\vp_i,\vp_j)\leq \gamma_n\]
then $\vp_0,\dots,\vp_{n+1}$ can be embedded into $\ell_2^n$: it is
\[ \gamma_n = 
\begin{cases}
\sqrt{1+\frac{2n+2}{n^2-2}} & \text{if $n$ is even},\\
\sqrt{1+\frac{2}{n-1}} & \text{if $n$ is odd}.
\end{cases}
\]
Their proof is geometrical.
The above analytical proof gives values that differ from $\gamma_n$ by $O(n^{-2})$ as $n\to\infty$.

\begin{corollary}[Sch\"utte \cite{Schutte}]\label{Schutte}
For any $\vp_1,\dots,\vp_{n+2}\in\ell_2^n$ we have
\[ \frac{\max\norm{\vp_i-\vp_j}_2}{\min_{i\neq j}\norm{\vp_i-\vp_j}_2} > 1+\frac{1}{n+2}.\]
\end{corollary}
Sch\"utte in fact determined the smallest $\delta_n>1$ such that
\[ \frac{\max\norm{\vp_i-\vp_j}_2}{\min_{i\neq j}\norm{\vp_i-\vp_j}_2} \geq 1+\delta_n\]
for all possible $\vp_1,\dots,\vp_{n+2}\in\ell_2^n$: it is $\delta_n=\gamma_{n+1}$.
This result was also discovered by Seidel \cite{Seidel}, and follows from the embedding result of Dekster and Wilker \cite{MR88k:51035}; see \cite{Barany} for a very simple proof.
\begin{corollary}\label{BMcor}
If $d(X,\ell_2^n)\leq 1+\frac{1}{n+2}$, then $e(X)\leq n+1$.
\end{corollary}
\begin{corollary}[Smyth \cite{Smyth}]\label{Smythcor}
If $\abs{p-2} < \frac{4(1+o(1))}{n\log n}$ then $e(\ell_p^n)=n+1$.
\end{corollary}
Compare this with Corollary~\ref{Swancor} in Section~\ref{lplowersection}.

\section{The Theorem of Brass and Dekster}\label{Brasssection}

\begin{theorem}[Brass \cite{MR2000i:52012} \& Dekster \cite{MR2001b:52001}]\label{Brassthm}
An $n$-dimensional Minkowski space contains an equilateral set of size $c(\log n)^{1/3}$ for some constant $c>0$ and $n$ sufficiently large.
\end{theorem}

Both Brass' and Dekster's proofs use Dvoretzky's theorem combined with a topological result.
We follow Brass' proof, which uses only the Brouwer fixed point theorem.

For the proof we need Dvoretzky's theorem.
\begin{Dvoretzky}
There exists a contant $c>0$ such that for any $\epsi>0$, any Minkowski space $X$ of sufficiently large dimension contains an $m$-dimensional subspace with Banach-Mazur distance at most $1+\epsi$ to $\ell_2^m$.
We may take $\dim X \geq e^{cm\epsi^{-2}}$ for some absolute constant $c>0$.
\end{Dvoretzky}
The estimate for $\dim X$, which is best possible (up to the value of $c$), was proven by Gordon \cite{MR87f:60058}.

\begin{proof}[Proof of Theorem~\ref{Brassthm}]
Denote the $n$-dimensional space by $X$.
Let $m=c_1(\log n)^{1/3}$ and $\epsi=1/(m+1)$.
A simple calculation gives that $e^{cm\epsi^{-2}} < n$ if we choose $c_1 = c_1(c)$ sufficiently small.
Dvoretzky's theorem gives that $X$ contains an $m$-dimensional subspace $Y$ with Banach-Mazur distance at most $1+\frac{1}{m+1}$ to $\ell_2^m$.
The proof is then finished by applying the following theorem.
\end{proof}

\begin{theorem}[Brass \cite{MR2000i:52012} \& Dekster \cite{MR2001b:52001}]
Let $X$ be an $n$-dimensional Minkowski space with Banach-Mazur distance $d(X,\ell_2^n)\leq 1+\frac{1}{n+1}$.
Then an equilateral set in $X$ of at most $n$ points can be extended to one of $n+1$ points.
In particular, $e(X)\geq n+1$.
\end{theorem}
\begin{proof}
By induction on nested subspaces of $X$ it is sufficient to prove that a $1$-equilateral set $\vp_1,\dots,\vp_n$ may be extended.
We may assume that a coordinate system has been chosen such that
\[ \norm{\vx}\leq\norm{\vx}_2\leq(1+\epsi)\norm{\vx},\]
where $\epsi=(n+1)^{-1}$.
Let $\vp_1,\dots,\vp_n$ be a $1$-equilateral set in $X$.
We want to find a $\vp\in X$ such that $\norm{\vp-\vp_i}=1$ for all $i=1,\dots,n$.

We have
\[ 1\leq\norm{\vp_i-\vp_j}_2\leq 1+\epsi\;\text{for all distinct }i,j.\]
By Theorem~\ref{simplexemb} $\vp_1,\dots,\vp_n$ are affinely independent.
For any $\rho_1,\dots,\rho_n\in[1,1+\epsi]$ there exists a point $\vx\in X$ with $\norm{\vx-\vp_i}_2=\rho_i$, $i=1,\dots,n$, by Theorem~\ref{simplexemb}.
Furthermore, by fixing a half space bounded by the hyperplane through $\vp_1,\dots,\vp_n$, there is a unique such $\vx$ in this half space (by simple linear algebra arguments as in the proof of of Theorem~\ref{CMsufficiency}; exercise).
We write $\vx=:\vf(\rho_1,\dots,\rho_n)$, and define $\mathbf{\fhi}:[0,\epsi]^n\to[0,\epsi]^n$ by sending $(\delta_1,\dots,\delta_n)\mapsto(y_1,\dots,y_n)$, where
\[ y_i=\delta_i+1-\norm{\vf(1+\delta_1,\dots,1+\delta_n)-\vp_i},\; i=1,\dots,n.\]
Since $1+\delta_i\leq 1+\epsi$, each $y_i$ is well-defined.
Also,
\begin{eqnarray*}
y_i &\leq& \delta_i +1 - (1+\epsi)^{-1}\norm{\vf(1+\delta_1,\dots,1+\delta_n)-\vp_i}_2\\
&=&\delta_i+1-(1+\epsi)^{-1}(1+\delta_i)\\
&=&(1-(1+\epsi)^{-1})(1+\delta_i)\\
&\leq&(1-(1+\epsi)^{-1})(1+\epsi) = \epsi,
\end{eqnarray*}
and
\begin{eqnarray*}
y_i &\geq&  \delta_i+1-\norm{\vf(1+\delta_1,\dots,1+\delta_n)-\vp_i}_2\\
&=& \delta_i+1-(1+\delta_i)=0,
\end{eqnarray*}
which shows that $\mathbf{\fhi}$ maps into $[0,\epsi]^n$.

By Brouwer's fixed point theorem, $\mathbf{\fhi}$ has a fixed point $(\delta_1,\dots,\delta_n)$, which gives
\[ \norm{\vf(1+\delta_1,\dots,1+\delta_n)-\vp_i} =1\;\text{ for all $i=1,\dots,n$.}\qedhere\]
\end{proof}
\begin{corollary}[Brass \cite{MR2000i:52012} \& Dekster \cite{MR2001b:52001}]
If $d(X,\ell_2^n)\leq 1+\frac{1}{n+2}$ then $e(X)=n+1$.
\end{corollary}
\begin{proof}
Combine the above theorem with Corollary~\ref{BMcor}.
\end{proof}

\section{The Linear Algebra Method}\label{linalgsection}
Linear algebra provides a very powerful counting tool.
The basic idea is that whenever one has a set $S$ of $m$ elements, and one wants to find an upper bound to $m$, one constructs a vector $\vv_s$ in some $N$-dimensional vector space for each element of $s\in S$.
If one can prove that the vectors constructed are linearly independent, then one immediately has the upper bound $m\leq N$.
If the vectors are not linearly independent, then one may try to prove that the square matrix $A^\Tr A$ has a high rank, where $A$ is the matrix with column vectors $\vv_s$, $s\in S$.

\subsection{Linear independence}\label{linindsection}
As a first example, one can prove that $e(\ell_2^n)\leq n+1$ as follows.
Let $S\subset\ell_2^n$ be $1$-equilateral, and define for each $\vs\in S$ a polynomial in $n$ variables
\[ P_{\vs}(\vx)=P_{\vs}(x_1,\dots,x_n):=1-\norm{\vx-\vs}_2^2.\]
Then it is easily seen that $\{1\}\cup\{P_{\vs}(\vx):\vs\in S\}$ is linearly independent in the vector space $\R[x_1,\dots,x_n]$ of polynomials in $x_1,\dots,x_n$ with real coefficients, and secondly that all $P_{\vs}(\vx)$ are in the linear span of $\{1,\sum_{i=1}^n x_i^2, x_1,\dots,x_n\}$.
It follows that $1+\card{S}\leq 2+n$, which is what was required.

In the same way one may prove that $e(\ell_p^n)\leq 1+(p-1)n$ if $p$ is an even integer (an observation due to Galvin; see \cite{Smyth}).
With some more work the following may be shown:
\begin{theorem}[Swanepoel \cite{Swanepoel-AdM}]
For each $n\geq 1$, $e(\ell_4^n)=n+1$.
\end{theorem}
\begin{proof}
The lower bound has already been shown in Section~\ref{kusner}.

Let $S$ be a $1$-equilateral set in $\ell_4^n$.
For each $\vs\in S$, let $P_{\vs}(x)=p_{\vs}(x_1,\dots,x_n)$ be the following polynomial:
\begin{eqnarray}
P_{\vs}(\vx) & := & 1-\norm{\vx-\vs}_4^4 \nonumber \\
& = & 1-\norm{\vs}_4^4-\sum_{i=1}^n x_i^4 + 4\sum_{i=1}^n s_i x_i^3 - 6\sum_{i=1}^n s_i^2 x_i^2 + 4\sum_{i=1}^n s_i^3 x_i.\label{zero}
\end{eqnarray}
Thus each $P_{\vs}$ is in the linear span of
\[\{1, \sum_{i=1}^n x_i^4\}\cup\{x_i^k: i=1,\dots,n; k=1,2,3\},\]
which is a subspace of $\R[x_1,\dots,x_n]$ of dimension $2+3n$.
Since $P_{\vs}(\vs)=1$ and $P_{\vs}(\vs')=0$ for all distinct $\vs,\vs'\in S$ we have that  $\{P_{\vs}:\vs\in S\}$ is linearly independent.
(So already $|S|\leq 2+3n$.)

We now prove that
\[ \{P_{\vs}: \vs\in S\}\cup\{1\}\cup\{x_i^k:i=1,\dots,n; k=1,2\} \]
is still linearly independent, which will give $|S|+1+2n\leq 2+3n$, proving the theorem.

Let $\lambda_{\vs}\, (\vs\in S), \lambda, \lambda_i, \mu_i\, (i=1,\dots,d) \in\R$ satisfy
\begin{equation}\label{one} \sum_{\vs\in S}\lambda_{\vs} P_{\vs}(\vx) + \lambda 1 +\sum_{i=1}^n\lambda_i x_i + \sum_{i=1}^n\mu_i x_i^2 = 0 \qquad \forall x_1,\dots,x_n\in\R.
\end{equation}
Substitute \eqref{zero} into \eqref{one}:
\begin{eqnarray*}
& & (-\sum_{\vs\in S}\lambda_{\vs})\sum_{i=1}^n x_i^4 + 4\sum_{i=1}^n(\sum_{\vs\in S}\lambda_{\vs}s_i)x_i^3 + \sum_{i=1}^n(\mu_i-6\sum_{\vs\in S}\lambda_{\vs}s_i^2)x_i^2 \\ & & + \sum_{i=1}^n(\lambda_i+4\sum_{\vs\in S}s_i^3)x_i +\sum_{\vs\in S}\lambda_{\vs} - \sum_{\vs\in S}\lambda_{\vs}\norm{\vs}_4^4+\lambda = 0.
\end{eqnarray*}
Since the above is a polynomial in $\R[x_1,\dots,x_n]$ which is identically $0$, the coefficients are also $0$:
\begin{eqnarray}
& & \sum_{\vs\in S}\lambda_{\vs} = 0, \label{two}\\
& & \sum_{\vs\in S}\lambda_{\vs}s_i=0 \qquad \forall i=1,\dots,n, \label{three}\\
& & \mu_i-6\sum_{\vs\in S}\lambda_{\vs}s_i^2 = 0 \qquad \forall i=1,\dots,n, \label{four}\\
& & \lambda_i+4\sum_{\vs\in S}\lambda_{\vs}s_i^3=0 \qquad \forall i=1,\dots,d, \label{five}\\
& & \lambda+\sum_{\vs\in S}\lambda_{\vs}-\sum_{\vs\in S}\lambda_{\vs}\norm{\vs}_4^4=0.\label{six}
\end{eqnarray}
By substituting $\vx=\vs\in S$ into \eqref{one} we obtain
\begin{equation}\label{seven}
\lambda_{\vs}+\lambda+\sum_{i=1}^n\lambda_i s_i+\sum_{i=1}^{n}\mu_i s_i^2 = 0\qquad \forall \vs\in S.
\end{equation}
Now multiply \eqref{seven} by $\lambda_{\vs}$ and sum over all $\vs\in S$:
\[ \sum_{\vs\in S}\lambda_{\vs}^2 + \lambda\sum_{\vs\in S}\lambda_{\vs} +\sum_{i=1}^n\lambda_i(\sum_{\vs\in S}\lambda_{\vs}s_i) + \sum_{i=1}^{n}\mu_i(\sum_{\vs\in S}\lambda_{\vs}s_i^2) = 0.\]
Then use \eqref{two}--\eqref{four} to simplify this expression as follows:
\[ \sum_{\vs\in S}\lambda_{\vs}^2 + \frac16\sum_{i=1}^n\mu_i^2 = 0.\]
It follows that all $\lambda_{\vs}=0$ and all $\mu_i=0$.
From \eqref{five} it follows that all $\lambda_i=0$ and from \eqref{six} that $\lambda=0$.

Thus we have linear independence.
\end{proof}

In the same way the following may be proved:
\begin{theorem}[Swanepoel \cite{Swanepoel-AdM}]\label{evenp}
For $p$ an even integer and $n\geq 1$ we have
\[ e(\ell_p^n) \leq \begin{cases} (\tfrac{p}{2}-1)n+1 & \text{if } p\equiv 0 \pmod{4},\\
\tfrac{p}{2}n+1 & \text{if } p\equiv 2 \pmod{4}. \end{cases} \]
\end{theorem}

\subsection{Rank arguments}\label{ranksection}
The following lemma is very simple, yet extremely powerful.
\begin{ranklemma}
Let $A=\begin{bmatrix}a_{ij}\end{bmatrix}$ be a real symmetric $n\times n$ non-zero matrix.
Then \[\rank A\geq \frac{(\sum_{i=1}^na_{ii})^2}{\sum_{i,j=1}^n a_{ij}^2}.\]
\end{ranklemma}
\begin{proof}
Let $r=\rank A$, and let $\lambda_1,\dots,\lambda_r\in\R$ be all non-zero eigenvalues of $A$.
Then $\sum_i a_{ii} = \trace(A)=\sum_{i=1}^r\lambda_i$.
Also, $\sum_{i,j}a_{ij}^2 = \trace(A^2)=\sum_{i=1}^r\lambda_i^2$, since $A^2$ has non-zero eigenvalues $\lambda_i^2$.
Therefore,
\[ \frac{(\sum_{i=1}^na_{ii})^2}{\sum_{i,j=1}^n a_{ij}^2}
=\frac{(\sum_{i=1}^r\lambda_i)^2}{\sum_{i,j=1}^r \lambda_i^2} \leq r\]
by the Cauchy-Schwartz inequality.
\end{proof}

\begin{corollary}\label{rankcor}
Let $A=\begin{bmatrix}a_{ij}\end{bmatrix}$ be a real symmetric $n\times n$ matrix with $a_{ii}=1$ for all $i$ and $\abs{a_{ij}}\leq\epsi$ for all $i\neq j$.
Then \[ \rank A \geq \frac{n}{1+(n-1)\epsi^2}.\]
In particular, if $\epsi=n^{-1/2}$ then $\rank A\geq n/2$.
\end{corollary}

Also, if $\epsi<1/(n-1)$ then the corollary gives that $A$ is invertible, which also follows from the fact that $A$ is then diagonally dominated.

The rank lemma together with a simple rounding argument gives an upper bound for $e(\ell_p^n)$, $1\leq p<\infty$, that is polynomial in $p$ and $n$.
This is due to Alon (Smyth, personal communication).
\begin{theorem}\label{firstlp}
For some constant $c>0$ we have that for $1\leq p<\infty$ and $n\geq 1$, $e(\ell_p^n)\leq cp^2n^{2+2/p}$.
\end{theorem}

\begin{proof}
Let $\vp_1,\dots,\vp_m$ be $1$-equilateral in $\ell_p^n$.
Since $\abs{p_i^{(k)}-p_j^{(k)}}\leq 1$ for all $i\neq j$ and all $k=1,\dots,m$, we may assume after a translation that all $p_i^{(k)}\in[0,1]$.
Choose an integer $N$ such that
\[ (e-1)pn^{1/p}\sqrt{m}\leq N < N+1 \leq epn^{1/p}\sqrt{m}.\]
(The interval $pn^{1/p}\sqrt{m}\geq 2$ since $m\geq 4$ without loss of generality.)
We now round each $\vp_i$ to a point in the lattice $\{0,\frac{1}{N}, \frac{2}{N},\dots,\frac{N}{N}\}^n$ as follows.
Let $q_i^{(k)}$ be the integer multiple of $1/N$ nearest to $p_i^{(k)}$, say $q_i^{(k)}=\frac{d(i,k)-1}{N}$ where $1\leq d(i,k)\leq N+1$ (if there is a tie, choose arbitrarily).
Let $\vq_i=(q_i^{(1)},\dots,q_i^{(n)})$, and let $Q$ be the $m\times m$ matrix $\begin{bmatrix}1-\norm{\vq_i-\vq_j}_p^p\end{bmatrix}_{i,j}^m$.
Thus $Q$ is an approximation of the identity $I_m$, and we now estimate how close.
The diagonal contains only $1$'s.
For each $k$, $\abs{p_i^{(k)}-q_i^{(k)}}\leq 1/2N$, and
\begin{eqnarray*}
\norm{\vq_i-\vq_j}_p &\leq& \norm{\vp_i-\vp_j}_p + \norm{\vp_i-\vq_i}_p+\norm{\vp_j-\vq_j}_p\\
&\leq& 1+n^{1/p}/N,
\end{eqnarray*}
and similarly,
$\norm{\vq_i-\vq_j}_p\geq 1-n^{1/p}/N$.
Thus
\[\norm{\vq_i-\vq_j}_p^p\leq\left(1+\frac{n^{1/p}}{N}\right)^p\leq\left(1+\frac{1}{(e-1)p\sqrt{m}}\right)^p < e^{\frac{1}{(e-1)\sqrt{m}}} < 1+\frac{1}{\sqrt{m}},\]
since $e^x<1+(e-1)x$ for $0<x<1$.
It follows that $1-\norm{\vq_i-\vq_j}_p^p > -1/\sqrt{m}$.
Since $(1+x)^p+(1-x)^p>2$ for $0<x<1$,
\[ 1-\norm{\vq_i-\vq_j}_p^p\leq 1-\left(1-\frac{n^{1/p}}{N}\right)^p<\left(1+\frac{n^{1/p}}{N}\right)^p-1< \frac{1}{\sqrt{m}}.\]
Thus by Corollary~\ref{rankcor}, $\rank Q\geq m/2$.
Note that $Q=\sum_{k=1}^n Q_k$, where
\[Q_k=\begin{bmatrix}\frac{1}{n}-\abs{q_i^{(k)}-q_j^{(k)}}^p\end{bmatrix}_{i,j=1}^m.\]
Define the $(N+1)\times m$ matrix
\[A_k=\begin{bmatrix} \ve_{d(1,k)} & \ve_{d(2,k)} & \cdots & \ve_{d(m,k)}\end{bmatrix},\]
 where $\ve_1,\dots,\ve_{N+1}$ is the standard basis of $\R^{N+1}$ (recall $q_i^{(k)} = \frac{d(i,k)-1}{N}$).
Secondly, define
\[B=\begin{bmatrix}\frac{1}{n}-\abs{\frac{i-1}{N}-\frac{j-1}{N}}^p\end{bmatrix}_{i,j=1}^{N+1}.\]
Then it is easily seen that $Q_k = A_k^{\Tr} B A_k$.
It follows that $\rank Q_k\leq N+1$, hence
\[\rank Q\leq\sum_{k=1}^n\rank Q_k\leq n(N+1).\]

Putting the upper and lower bounds for $\rank Q$ together we obtain
\[ m/2 < n(N+1) \leq epn^{1+1/p}\sqrt{m},\]
and $m<4e^2p^2n^{2+2/p}$.
\end{proof}

\section{Approximation Theory: Smyth's approach}
In the proof of Theorem~\ref{firstlp} we essentially approximated the function $f_p(t)=\abs{t}^p$, $t\in[-1,1]$ uniformly by step functions (look at the definition of the matrix $B$ in the proof).
For $p>1$ the function $f_p$ is differentiable $\lfloor p\rfloor$ times, so one would expect that there are better uniform approximations to $f_p$ using polynomials or splines instead of step functions.
This is indeed the case, and can be used to improve the bound of Theorem~\ref{firstlp}.
Smyth \cite{Smyth} used the theorems of Jackson from approximation theory to find the first non-trivial upper bounds for $e(\ell_p^n)$.
Alon and Pudl\'ak \cite{MR1995795} improved his bounds using the rank lemma.
We state the theorem of Jackson that is needed (called Jackson V) and then prove their result.
The proofs of the Jackson theorems may be found in many texts on approximation theory, e.g.\ \cite{MR99f:41001}.

\begin{JacksonV}
There is an absolute constant $c>0$ such that for any $f\in C[-1,1]$
\begin{enumerate}
\item there exists a polynomial $P$ of degree at most $n$ such that
\[ \norm{f-P}_\infty \leq c\omega(f,1/n),\]
where $\omega(f,\delta):=\sup\{\abs{f(x)-f(y)}: \abs{x-y}\leq\delta\}$ is the \emph{modulus of continuity} of $f$,
\item and if $f$ is $k$ times differentiable and $n\geq k$, then there exists a polynomial $P$ of degree at most $n$ such that
\[ \norm{f-P}_\infty \leq \frac{c^k}{(n+1)n(n-1)\dots(n-k+2)}\omega(f^{(k)},1/n).\]
\end{enumerate}
\end{JacksonV}

\begin{lemma}\label{plemma}
For each $1\leq p<\infty$ there exists a constant $c_p>0$ such that for any $d\geq 1$ there is a polynomial $P$ of degree at most $d$ such that $P(0)=0$ and
\[ \abs{\abs{t}^p-P(t)}\leq\frac{c_p}{d^p}\quad\text{for all $t\in[-1,1]$.}\]
We have $c_p < (cp)^p$ for some universal $c>0$.
\end{lemma}
\begin{proof}
By choosing $c_p$ sufficiently large we may assume that $d>2p$.
Then $f_p(t):=\abs{t}^p$ satisfies $f\in C^{(k)}[-1,1]$ with $k=\lceil p\rceil -1$, and
\[ f^{(\lceil p\rceil -1)}(t)=\sgn^{\lceil p\rceil -1}(t)p(p-1)\dots(p-\lceil p\rceil+2)\abs{t}^{p+1-\lceil p\rceil}.\]
Thus \[\omega(f^{(\lceil p\rceil -1)}, \delta)=p(p-1)\dots(p-\lceil p\rceil+2)\delta^{p+1-\lceil p\rceil}.\]
From Jackson V we obtain a polynomial $P$ of degree at most $d$ with
\begin{eqnarray*}
\norm{f_p-P}_\infty &\leq& \frac{c^{\lceil p\rceil -1}}{(d+1)d(d-1)\dots(d-\lceil p\rceil +3)}p(p-1)\dots(p-\lceil p\rceil+2)(1/d)^{p+1-\lceil p\rceil}\\
&<& \frac{c^{\lceil p\rceil -1}p(p-1)\dots(p-\lceil p\rceil+2)}{(\frac12 d)^{\lceil p\rceil -1}}(1/d)^{p+1-\lceil p\rceil}\\
&=& \frac{c}{d^p}.
\end{eqnarray*}
If we subtract the constant term from $P$ we obtain a polynomial $Q$ with $Q(0)=0$ and $\norm{f_p-Q}_\infty\leq 2c/d^p$.
\end{proof}

Smyth \cite{Smyth} approached Kusner's problem with the idea of approximating the $p$-norm with polynomials using Jackson's theorems.
He obtained an upper bound $c_p n^{(p+1)/(p-1)}$ for $1<p<\infty$ using linear independence.
The next theorem was proved by Alon and Pudl\'ak by combining Smyth's approach with the rank lemma.

\begin{theorem}[Smyth \cite{Smyth}, Alon \& Pudl\'ak \cite{MR1995795}]
For $1\leq p<\infty$, $e(\ell_p^n)\leq c'_p n^{(2p+2)/(2p-1)}$.
We may take $c'_p=cp$, with $c>0$ and absolute constant.
\end{theorem}

\begin{proof}
Let $\vp_1,\dots,\vp_m\in\ell_p^n$ be $1$-equilateral.
Let $c=\max(c_p,(2^{1/p}-1)^{-p})$, where $c_p$ is the constant from Lemma~\ref{plemma}.
We fix an integer $d$ such that $cn\sqrt{m} < d^p < 2cn\sqrt{m}$ (possible since $c\geq (2^{1/p}-1)^{-p}$).
Let $P$ be the polynomial from Lemma~\ref{plemma}.
Define the $m\times m$ matrix $A=\begin{bmatrix}a_{ij}\end{bmatrix}$ by $a_{ij}=1-\sum_{k=1}^nP(p_i^{(k)}-p_j^{(k)})$.
Since $\sum_{k=1}^nP(p_i^{(k)}-p_j^{(k)})$ is an approximation of $\norm{\vp_i-\vp_j}_p^p$, the matrix $A$ is an approximation of the identity $I_m$.
Again we estimate how close.
Firstly, since $P(0)=0$, all $a_{ii}=1$.
Secondly for any $i\neq j$,
\begin{eqnarray*}
\abs{a_{ij}} &=& \bigabs{\norm{\vp_i-\vp_j}_p^p - \sum_{k=1}^nP(p_i^{(k)}-p_j^{(k)})}\\
&\leq& \sum_{k=1}^n\bigabs{\abs{\vp_i^{(k)}-\vp_j^{(k)}}^p - P(\vp_i^{(k)}-\vp_j^{(k)})}\\
&\leq& n\frac{c_p}{d^p} < \frac{1}{\sqrt{m}},
\end{eqnarray*}
by choice of $d$.
By Corollary~\ref{rankcor} we have $\rank A \geq m/2$.

We now find an upper bound for the rank of $A$.
For each $i=1,2,\dots,m$, define the polynomial in $n$ variables
\[ P_i(x_1,\dots,x_n) := 1- \sum_{k=1}^n P(p_i^{(k)}-x_i).\]
Thus $a_{ij}=P_i(\vp_j)$.
Each $P_i$ is in the linear span of
\[ \mathcal{P} = \{1,\sum_{k=1}^n x_i^d, x_1,\dots,x_n,x_1^2,\dots,x_n^2,\dots,x_1^{d-1},\dots,x_n^{d-1}\}.\]
Thus the $m$ polynomials $P_i$ lie in a $(2+(d-1)n)$-dimensional subspace of polynomials.
Then the $i$th row vector of $A$, $\begin{bmatrix}P_i(\vp_1), \dots, P_i(\vp_m)\end{bmatrix}$, is in the linear span of
\[\{ (f(\vp_1),\dots,f(\vp_m)): f\in\mathcal{P}\},\]
which is a subspace of $\R^m$ of dimension at most $\card{\mathcal{P}}=2+(d-1)n\leq dn$.
It follows that $\rank A \leq dn$.

Putting the upper and lower bound for $\rank A$ together, we find $m/2 \leq dn$.
Since $d^p < 2cn\sqrt{m}$, we obtain $m < c'_p n^{(2p+2)/(2p-1)}$.
\end{proof}

In the next section we find the best known upper bounds for $e(\ell_1^n)$, due to Alon and Pudl\'ak.

\section{The best known upper bound for $e(\ell_1^n)$}
Alon and Pudl\'ak \cite{MR1995795} proved that $e(\ell_p^n)\leq c_pn\log n$ if $p$ is an odd integer.
This matches the lower bound of $2n$ apart from the $\log n$ factor and the constant that depends on $p$.
We here present their proof for the case $p=1$.
The proof for other odd $p$ is simple once this case is understood, see \cite{MR1995795} for the detail.

\begin{theorem}[Alon and Pudl\'ak \cite{MR1995795}]
$e(\ell_1^n)\leq cn\log n$.
\end{theorem}
\begin{proof}
Let $\vp_1,\dots,\vp_{m+1}$ be $1$-equilateral in $\ell_1^n$.
After a translation we may assume $\vp_{m+1}=\vo$.
Thus it is sufficient to prove $m+1\leq cn\log n$ given that
\begin{eqnarray}
\sum_{k=1}^n\abs{p_i^{(k)}}&=&1\quad\forall i=1,\dots,m,\label{star1}\\
\text{and } \sum_{k=1}^n\abs{p_i^{(k)}}&=&1\quad\forall 1\leq i\neq j\leq m.\label{star2}
\end{eqnarray}
We first show that by doubling the dimension we may assume that all $p_i^{(k)}\geq 0$.
We replace each $p_i^{(k)}$ by
\[ \begin{cases} (p_i^{(k)},0) & \text{if $p_i^{(k)}\geq 0$,}\\ (0,-p_i^{(k)}) & \text{if $p_i^{(k)}\leq 0$}.\end{cases}\]
Now the points $\vp_1,\dots,\vp_m$ are in $\ell_1^{2n}$, they are still $1$-equilateral, $\norm{\vp_i}_1=1$, and now all $p_i^{(k)}\geq 0$.
Thus we may assume together with \eqref{star1} and \eqref{star2} that $p_i^{(k)}\geq 0$ holds.
(Once we have proven $m\leq cn\log n$ with this assumption, then $m\leq 1+c(2n)\log(2n)<3cn\log n$ in the general case.)

We now express $\abs{p_i^{(k)}-p_j^{(k)}}$ in terms of $\min(p_i^{(k)},p_j^{(k)})$.
Since
\[ \abs{a-b}=a+b-2\min(a,b)\quad\text{for $a,b\in\R$,}\]
and $\sum_{k=1}^n p_i^{(k)}=1$, we obtain
\[ I_m = \begin{bmatrix}1-\norm{\vp_i-\vp_j}\end{bmatrix}_{i,j=1}^m = \begin{bmatrix}\displaystyle -1+2\sum_{k=1}^n\min\big(p_i^{(k)},p_j^{(k)}\big)\end{bmatrix}_{i,j=1}^m.\]
We now want to approximate $\min(a,b)$ by an inner product $\ipr{\bva}{\bvb}$ where $\bva,\bvb\in\R^N$.
Then $I_m$ will be approximated by
\[ A=\begin{bmatrix}-1+2\ipr{\big(\bvp_i^{(1)},\dots,\bvp_i^{(n)}\big)}{\big(\bvp_j^{(1)},\dots,\bvp_j^{(n)}\big)}\end{bmatrix}_{i,j=1}^m.\]
Then $\rank A\leq 1+Nn$, since it is a linear combination of $J_m$ of rank $1$ and a Gram matrix of vectors in $\R^{Nn}$.
The approximation will be good enough so that the rank lemma will give $\rank A>c_1m$.
Thus $c_1m<Nn+1$, and we will have to make sure that $N<c_2\log n$.
What will in fact happen is that $N$ will be different for each coordinate $k=1,\dots,n$, say $N_k$, and we'll have to make sure $\sum_{k=1}^n N_k < c_2n\log n$.

The rough idea for approximating $\min(a,b)$ by an inner product is the following.
For each $x\in[0,1]$ define
\[\bvx = (\underbrace{\tfrac{1}{\sqrt{N}},\dots,\tfrac{1}{\sqrt{N}},}_{\lfloor x N\rfloor\text{ times}} 0,\dots,0)\in\R^N.\]
Then for $0\leq a\leq b\leq 1$,
\[\ipr{\bva}{\bvb}=\tfrac{1}{N}\lfloor aN\rfloor = a-\epsi=\min(a,b)-\epsi,\]
where $0\leq\epsi\leq1/N$.
Let us see how far this brings us already.
Since
\[ \bigabs{\sum_{k=1}^n\min\big(p_i^{(k)},p_j^{(k)}\big) - \ipr{\big(\bvp_i^{(1)},\dots,\bvp_i^{(n)}\big)}{\big(\bvp_j^{(1)},\dots,\bvp_j^{(n)}\big)}} \leq \frac{n}{N},\]
the error by which $A$ approximates $I_m$ is
\[ a_{ij}-\delta_{ij} \leq \frac{2n}{N}=\frac{1}{\sqrt{m}}\]
if we choose $N=2n\sqrt{m}$.
By the rank lemma, $\rank A\geq (1+o(1))m/2$.
Also, $\rank A\leq Nn+1$.
Putting the bounds together, we obtain $m<cn^4$, already obtained in Theorem~\ref{firstlp}.

We now fine-tune this idea.
For any partition $\CP=\{0=u_0<u_1<\dots<u_{N-1}<u_N=1\}$ of $[0,1]$ into $N$ intervals, we define for any $x\in[0,1]$ the vector $\vx_{\CP}\in\R^N$, where
\[ x_{\CP}^{(j)}=\begin{cases} \sqrt{u_j-u_{j-1}}, & j<t, \\ \frac{x-u_{t-1}}{\sqrt{u_t-u_{t-1}}}, & j=t,\\ 0, & j>t, \end{cases}\]
where $t$ is such that $x\in (u_{t-1},u_t)$.

\subsubsection*{Remarks.}
\begin{enumerate}
\item We'll choose a different partition for each coordinate $k=1,\dots,n$.
Since both $\sum_k\min(p_i^{(k)},p_j^{(k)})$ and $\ipr{(\bvp_i^{(1)},\dots,\bvp_i^{(n)})}{(\bvp_j^{(1)},\dots,\bvp_j^{(n)})}$ are sums over $k=1,\dots,n$, the error in each coordinate will also sum up, so we first do the analysis in a single coordinate.
\item We'll choose the partitions so that a coordinate will never hit an endpoint of any $(u_{j-1},u_j)$.
\end{enumerate}

With the above definition we find for $0\leq a\leq b\leq 1$ that
\begin{eqnarray*}
\ipr{\va_{\CP}}{\vb_{\CP}} &=& (\sqrt{u_1-0})^2 + (\sqrt{u_2-u_1})^2 + (\sqrt{u_{t-1}-u_{t-2}})^2 + a_{\CP}^{(t)}b_{\CP}^{(t)},\quad\text{where $a\in(u_{t-1},u_t)$,}\\
&=&\begin{cases} u_{t-1} + \frac{a-u_{t-1}}{\sqrt{u_t-u_{t-1}}}\sqrt{u_t-u_{t-1}} & \text{if $b\notin(u_{t-1},u_t)$,}\\ u_{t-1} + \frac{a-u_{t-1}}{\sqrt{u_t-u_{t-1}}}\frac{b-u_{t-1}}{\sqrt{u_t-u_{t-1}}} & \text{if $b\in(u_{t-1},u_t)$,}\end{cases}\\
&=&\begin{cases} \min(a,b) & \text{if $a$ and $b$ are in different intervals of $\CP$,}\\ u_{t-1}+\frac{(a-u_{t-1})(b-u_{t-1})}{u_t-u_{t-1}} & \text{if $a,b\in(u_{t-1},u_t)$.}\end{cases}
\end{eqnarray*}
Thus there is no error if $a$ and $b$ are in different intervals.
However, when $a$ and $b$ are in the same interval of $\CP$, the error is bad, and we want to get rid of it.
We do this by adding another $N$ coordinates:
Define
\[ \hvx_{\CP} = \frac{x-u_{t-1}}{\sqrt{u_t-u_{t-1}}}\ve_t\in\R^N,\quad\text{with $t$ such that $x\in(u_{t-1},u_t)$.}\]
Then
\[\ipr{\hva_{\CP}}{\hvb_{\CP}}= \begin{cases} 0 & \text{if $a$ and $b$ are in different intervals of $\CP$,}\\ \frac{(a-u_{t-1})(b-u_{t-1})}{u_t-u_{t-1}} & \text{if $a,b\in(u_{t-1},u_t)$.}\end{cases}\]
Since we want to subtract this, we let the inner product on $\R^N\oplus\R^N$ be
\[ \ipr{(\vx_1,\vx_2)}{(\vy_1,\vy_2)} = \ipr{\vx_1}{\vy_1}-\ipr{\vx_2}{\vy_2}.\]
We then obtain for $0\leq a\leq b\leq 1$ that
\[\ipr{(\va_{\CP},\hva_{\CP})}{(\vb_{\CP},\hvb_{\CP})}= \begin{cases} \min(a,b) & \text{if $a$ and $b$ are in different intervals of $\CP$,}\\ u_{t-1} & \text{if $a,b\in(u_{t-1},u_t)$.}\end{cases}\]
Thus if $a$ and $b$ are in the same interval, the approximation is a systematic rounding down.
We now add a third vector in $\R^N$ to improve the approximation in this case.
We do this by \emph{randomized rounding}:
For each interval $(u_{t-1},u_t)$ we choose randomly and uniformly $\tau\in(u_{t-1},u_t)$ (a \emph{threshold}), and for any $x\in(u_{t-1},u_t)$ we define
\[ \vvx_{\CP} = \begin{cases} 0 & \text{if $x\leq\tau$,}\\ \sqrt{u_t-u_{t-1}}\ve_t & \text{if $x>\tau$.}\end{cases} \]
Note that $\vvx_{\CP}$ is a random variable.
Thus if $a$ and $b$ are in different intervals then $\ipr{\vva_{\CP}}{\vvb_{\CP}}=0$, and if $a$ and $b$ are in the same interval $(u_{t-1},u_t)$, then (as with $\hvx_{\CP}$)
\[ \ipr{\vva_{\CP}}{\vvb_{\CP}} =\begin{cases} 0 & \text{if $\min(a,b)\leq\tau$,}\\ u_t-u_{t-1} & \text{if $\min(a,b)>\tau$.}\end{cases} \]
We now let the inner product on $\R^N\oplus\R^N\oplus\R^N$ be
\[ \ipr{(\vx_1,\vx_2,\vx_3)}{(\vy_1,\vy_2,\vy_3)} = \ipr{\vx_1}{\vy_1}-\ipr{\vx_2}{\vy_2}+\ipr{\vx_3}{\vy_3},\]
and for each $x\in[0,1]$, let
\[ \bvx_{\CP} = (\vx_{\CP},\hvx_{\CP},\vvx_{\CP})\in\R^{3N}.\]
Then for any $a,b\in[0,1]$,
\[ \ipr{\bva_{\CP}}{\bvb_{\CP}} = \begin{cases} \min(a,b) & \text{if $a,b$ are in different intervals of $\CP$,}\\ \begin{cases} u_{t-1}, & \text{if }\min(a,b)\leq\tau\\ u_t, &\text{if }\min(a,b)>\tau\end{cases}, & \text{if $a,b\in(u_{t-1},u_t)$,}\end{cases}\]
and the error $X = \min(a,b)-\ipr{\bva_{\CP}}{\bvb_{\CP}}$ (depending on $\CP, a, b$) satisfies
\[ X = \begin{cases} 0 & \text{if $a,b$ in different intervals,}\\ a-u_{t-1} & \text{if $\min(a,b)\leq\tau,\; a,b\in(u_{t-1},u_t)$,}\\ a-u_t & \text{if $\min(a,b)>\tau,\; a,b\in(u_{t-1},u_t)$.}\end{cases}\]
Thus \[\abs{X}\leq u_t-u_{t-1} = \length{u_{t-1},u_t},\] the length of the interval.
We now calculate the expected value and variance of $X$.
Without loss of generality $a\leq b$.
Then easily
\begin{eqnarray*}
E(X)&=&\prob(a\leq\tau)(a-u_{t-1})+\prob(a>\tau)(a-u_{t-1})\\
&=& \frac{u_t-a}{u_t-u_{t-1}}(a-u_{t-1}) + \frac{a-u_{t-1}}{u_t-u_{t-1}}(a-u_t)\\
&=& 0,
\end{eqnarray*}
and
\begin{eqnarray*}
E(X^2)&=&\prob(a\leq\tau)(a-u_{t-1})^2+\prob(a>\tau)(a-u_{t-1})^2\\
&=& \frac{u_t-a}{u_t-u_{t-1}}(a-u_{t-1})^2 + \frac{a-u_{t-1}}{u_t-u_{t-1}}(a-u_t)^2\\
&=& \frac{(a-u_{t-1})(u_t-a)}{u_t-u_{t-1}}(a-u_{t-1}-a+u_t)\\
&=& (a-u_{t-1})(u_t-a)\\
&\leq& \frac14(u_t-u_{t-1})^2, \text{ since } u_{t-1}<a<u_t.
\end{eqnarray*}

For each coordinate $k=1,\dots,n$, we'll (soon) choose a different partition $\CP_k$ of $[0,1]$ into $N_k$ intervals, with an independent random threshold for each interval in each $\CP_k$.
For each $\vx=(x(1),\dots,x(n))$ we let
\[\bvx=(\bvect{x(1)}_{\CP_1},\dots,\bvect{x(n)}_{\CP_n})\in\R^{3N_1}\oplus\dots\oplus\R^{3N_n}.\]
Thus for each $\vp_i$ we now have a $\bvp_i=(\bvp_i^{(1)},\dots,\bvp_i^{(n)})\in\R^{3\sum_{k=1}^n N_k}$ (where in each coordinate a different partition is used, but not denoted anymore).
Then we approximate the identity
\[ I_m=\begin{bmatrix}\displaystyle -1+2\sum_{k=1}^n\min\big(p_i^{(k)},p_j^{(k)}\big)\end{bmatrix}_{i,j=1}^m\]
by
\[ A=\begin{bmatrix}-1+2\ipr{\bvp_i}{\bvp_j}\end{bmatrix}_{i,j=1}^m = \begin{bmatrix}\displaystyle -1+2\sum_{k=1}^n\ipr{\bvp_i^{(k)}}{\bvp_j^{(k)}}\end{bmatrix}_{i,j=1}^m\]
and we let $X_{i,j,k}=\min\big(p_i^{(k)},p_j^{(k)}\big)-\ipr{\bvp_i^{(k)}}{\bvp_j^{(k)}}$.
Note that $\{X_{i,j,k} : k=1,\dots,n\}$ are independent random variables for fixed $i,j$, since $X_{i,j,k}\equiv 0$ if $p_i^{(k)}$ and $p_j^{(k)}$ are in different intervals of $\CP_k$, or depends on a single threshold, with these thresholds independent as $k=1,\dots,n$.

We want to apply the rank lemma to $A$, so we start estimating $(\sum_i a_{ii})^2$ and $\sum_{i,j}a_{ij}^2$.
First of all,
\begin{eqnarray*}
a_{ii}&=&-1+2\sum_{k=1}^n\ipr{\bvp_i^{(k)}}{\bvp_i^{(k)}}\\
&=& -1+2\sum_{k=1}^n\min\big(p_i^{(k)},p_i^{(k)}\big) - 2\sum_{k=1}^n X_{i,i,k}\\
&=& 1 - 2\sum_{k=1}^n X_{i,i,k}.
\end{eqnarray*}
Similarly, for $i\neq j$,
\[ a_{ij}=-2\sum_{k=1}^n X_{i,j,k}.\]
Since $\abs{X_{i,i,k}}\leq\length{I_k}$, where $I_k$ is the interval of $\CP_k$ containing $p_i^{(k)}$, we obtain
\[\sum_{k=1}^n X_{i,i,k}\leq\sum_{k=1}^n\length{I_k},\]
and
\begin{eqnarray}
\sum_{i=1}^m a_{ii} &=& m-2\sum_{i=1}^m\sum_{k=1}^n X_{k,i,i}\notag\\
&\geq & m-2\sum_{i=1}^m\sum_{k=1}^n\length{I_k}\notag\\
&=& m-2\sum_{I}\hit{I}\length{I}\label{apone},
\end{eqnarray}
wehre the last sum is over all intervals of all $\CP_k$, and $\hit{I}$ is the number of $p_i^{(k)}$ ``hitting'' $I$, i.e., the number of ordered pairs $(i,k)$ with $p_i^{(k)}\in I$.

Secondly, for $i\neq j$, $a_{ij}^2=4(\sum_{k=1}^n X_{i,j,k})^2$, 
hence $E(a_{ij}^2)=4\sum_{k=1}^n E(X_{i,j,k}^2)$ due to independence and $E(X_{i,j,k})=0$.
Similarly,
\begin{eqnarray*}
E(a_{ii}^2) &=& 1-4\sum_{k=1}^n E(X_{i,i,k}) + 4\sum_{k=1}^n E(X_{i,i,k}^2)\\
&=& 1+ 4\sum_{k=1}^n E(X_{i,i,k}^2).
\end{eqnarray*}
Thus
\begin{eqnarray*}
E(\sum_{i,j=1}^m a_{ij}^2) &=& \sum_{i=1}^m E(a_{ii}^2) + \sum_{\substack{i,j=1\\ i\neq j}}^m E(a_{ij}^2)\\
&=& m+4\sum_{i=1}^m\sum_{k=1}^n E(X_{i,i,k}^2) + 4\sum_{\substack{i,j=1\\ i\neq j}}^m\sum_{k=1}^n E(X_{i,j,k}^2)\\
&=& m + 4\sum_{i,j=1}^m\sum_{k=1}^n E(X_{i,j,k}^2)\\
&\leq& m + \sum_{i,j=1}^m\sum_{k=1}^n \length{I_{i,j,k}}^2, \quad\parbox[t]{7cm}{with $I_{i,j,k}\in\CP_k$ such that $p_i^{(k)},p_j^{(k)}\in I_{i,j,k}$ if it exists, otherwise take $\length{I_{i,j,k}}=0$,}\\
&=& m+\sum_I \hit{I}^2\length{I}^2,
\end{eqnarray*}
where the last sum is again over all intervals.
Because we have found an upper bound on the expected value of $\sum_{i,j=1}^m a_{ij}^2$, it follows that there exists a choice of thresholds such that
\begin{equation}\label{aptwo}
 \sum_{i,j=1}^m a_{ij}^2 \leq m+\sum_I \hit{I}^2\length{I}^2.
\end{equation}
Fix such a choice of thresholds.
(So at this stage we leave randomness behind.)
The rank lemma together with \eqref{apone} and \eqref{aptwo} gives
\[\rank A \geq \frac{(m-2\sum_I \hit{I}\length{I})^2}{m+\sum_I \hit{I}^2\length{I}^2}\]
if
\begin{equation}\label{aprequire}
m-2\sum_I \hit{I}\length{I}\geq 0.
\end{equation}
On the other hand, $\rank A \leq 1+\sum_{k=1}^n 3N_k = 1+3\sum_I 1$,
thus
\begin{equation}\label{apthree}
\frac{(m-2\sum_I \hit{I}\length{I})^2}{m+\sum_I \hit{I}^2\length{I}^2}\leq 1+3\sum_I 1.
\end{equation}

Finally, we fix the partition $\CP_k$ for each coordinate $k$.
We assume that for each $k=1,\dots,n$, all $p_i^{(k)}$ are distinct.
(This is not essential: the $p_i^{(k)}$ may all be perturbed by a sufficiently small amount so as to weaken the left-hand side of \eqref{apthree} by at most $\epsi>0$, for any $\epsi>0$.)

We now assume that both $m$ and $n$ are powers of $4$.
If we can prove $m\leq cn\log n$ in this case, then $m\leq 16cn\log n$ in general, since we may round $m$ down to a power of $4$ (thus dividing $m$ by at most $4$), and $n$ up to a power of $4$ (thus multiplying $n$ by at most $4$).
All the logarithms until the end of the proof will be base $2$ (and indicated as such).
We first divide $[0,\frac{1}{n}]$ into $\sqrt{m/n}$ equal parts (thus each of length $1/\sqrt{mn}$).
These are called \emph{base intervals}.
Without loss of generality no $p_i^{(k)}$ is an endpoint of a base interval, again by making an infinitesimal perturbation.

Then start at $t=1$ and let $t$ go down until
\[\hit{[t,1]}\length{[t,1]}\geq \frac{cn\binlog n}{m}\] or until $t=2^{-1/3}$, whichever comes first.
(The constant $c$ will be fixed later; $c=10000$ is sufficient.)
If the stopping point $t_1$ equals some $p_i^{(k)}$, we go down slightly more, without increasing $\hit{[t,1]}$.
Then start with a new interval at $t_1$ and go down until
\[\hit{[t,t_1)}\length{[t,t_1]}\geq \frac{cn\binlog n}{m}\]
 or until $t=2^{-2/3}$, whichever comes first.
In general we go down until
\[\hit{[t,t_{s-1})}\length{[t,t_{s-1})}\geq \frac{cn\binlog n}{m}\]
 or until $t=2^{-s/3}$, for each $s=1,\dots,3\binlog n$.

If $\hit{I}\length{I}\geq cn\binlog n/m$, we call $I$ a \emph{regular interval}, otherwise a \emph{singular interval}.

Thus for a singular interval $I$ we have
\[\hit{I}\length{I}< \frac{cn\binlog n}{m}, \quad\text{$I$ singular.}\]
For a regular interval, if the stopping point did not hit a $p_i^{(k)}$, then $\hit{I}\length{I}= cn\binlog n/m$.
If the stopping point did hit a $p_i^{(k)}$, then $\hit{I}\length{I} > cn\binlog n/m$, but
\[(\hit{I}-1)\length{I}\leq \frac{cn\binlog n}{m}.\]
Thus
\[ \hit{I}\length{I} \geq \frac{cn\binlog n}{m}, \quad\text{$I$ regular,}\]
and
\[ \hit{I}\length{I} \leq \frac{2cn\binlog n}{m}\quad\text{if $\hit{I}\geq 2$ for regular $I$.}\]
Later we'll have to take care of regular $I$ with $\hit{I}=1$ separately: for them 
\[ \sum_{\substack{\text{$I$ regular}\\ \hit{I}=1}}\length{I}^2\leq \sum_{\substack{\text{$I$ regular}\\ \hit{I}=1}}\length{I}  < \sum_I\length{I}=n.\]

There are at most $3n\binlog n$ singular intervals.
We now bound the number of regular intervals.
For each regular interval $I$, let $s_I$ be the $s$ for which $I\subseteq[2^{-s/3},2^{-(s-1)/3}]$.
Then
\begin{eqnarray}
\sum_{\text{regular $I$}} 2^{-s_I/3} \hit{I} &\leq & \sum_{\text{regular $I$}}\min(I)\hit{I}\notag\\
&<& \sum_I\sum_{\substack{i=1\\ p_i^{(k)}\in I}}^m\min(I),\quad\text{where $k$ is such that $I\in\CP_k$,}\notag\\
&<& \sum_I\sum_{\substack{i=1\\ p_i^{(k)}\in I}}^m p_i^{(k)}=\sum_{i=1}^m\sum_{k=1}^n p_i^{(k)}=m.\label{star}
\end{eqnarray}
We now bound $\sum_{\text{regular $I$}}1/(2^{-s_I/3}\hit{I})$.
For each $s=1,\dots,3\binlog n$,
\begin{eqnarray*}
\sum_{\substack{\text{regular $I$}\\ s_I=s}}\length{I} &\leq& (2^{-(s-1)/3}-2^{-s/3})n\quad\text{(since for each coordinate the $I$'s don't overlap)}\\
&=& (2^{1/3}-1)2^{2/3}n < \tfrac13 2^{-s/3}n.
\end{eqnarray*}
Since $\length{I}\hit{I}\geq cn\binlog n/m$ for regular $I$,
\[ \frac{1}{\hit{I}}\leq\frac{m\length{I}}{cn\binlog n}.\]
Thus
\begin{eqnarray*}
\sum_{\substack{\text{regular $I$}\\ s_I=s}}\frac{1}{2^{-s/3}\hit{I}} &\leq& \frac{m}{2^{-s/3}cn\binlog n}\sum_{\substack{\text{regular $I$}\\ s_I=s}}\length{I}\\
&<& \frac{m}{3c\binlog n}.
\end{eqnarray*}
Thus
\begin{equation}\label{starstar}
\sum_{\text{regular $I$}}\frac{1}{2^{-s_I/3}\hit{I}} < \sum_{s=1}^{3\binlog n}\frac{m}{3c\binlog n}=\frac{m}{c}.
\end{equation}
By Cauchy-Schwartz,
\begin{eqnarray*}
\left(\sum_{\text{regular $I$}} 1 \right)^2 &=& \left(\sum_{\text{regular $I$}}\sqrt{2^{-s_I/3}\hit{I}}\frac{1}{\sqrt{2^{-s_I/3}\hit{I}}}\right)^2\\
&\leq& \left(\sum_{\text{regular $I$}} 2^{-s_I/3} \hit{I}\right)\left(\sum_{\text{regular $I$}}\frac{1}{2^{-s_I/3}\hit{I}}\right)\\
&<& m\frac{m}{c}\quad\text{by \eqref{star} and \eqref{starstar},}
\end{eqnarray*}
and the number of regular intervals is less than $m/\sqrt{c}$.
Thus the total number of non-base intervals is at most $3n\binlog n+m/\sqrt{c}$ (and recall that the number of base intervals is $\sqrt{mn}$).
We are now in a position to estimate the various quantities in \eqref{apthree}.
We want to show that these estimates imply $m\leq cn\binlog n$, so we assume that $m > cn\binlog n$ and aim for a contradiction if $c$ is sufficiently large ($10000$ will do).
First of all,
\begin{eqnarray}
1+3\sum_I 1 &<& 1+3\left(\sqrt{mn} + 3n\binlog n + \frac{m}{\sqrt{c}}\right)\notag\\
&<& 1 + \frac{3m}{\sqrt{c\binlog n}} + \frac{9m}{c} +\frac{3m}{\sqrt{c}}\notag\\
&<& \frac{16m}{\sqrt{c}}\label{eqa}
\end{eqnarray}
to be sure.
Secondly,
\begin{eqnarray}
m+\sum_I \hit{I}^2\length{I}^2 &=& m + \sum_{\text{base $I$}}\hit{I}^2\frac{1}{mn} + \sum_{\substack{\text{regular $I$}\\ \hit{I}=1}}\length{I}^2 \notag\\
&& + \sum_{\substack{\text{non-base $I$}\\ \text{$h_I\geq 2$ if regular}}}(\hit{I}\length{I})^2\notag\\
&<& m + nm^2\frac{1}{mn} + n + \left(3n\binlog n +\frac{m}{\sqrt{c}}\right)\left(\frac{2cn\binlog n}{m}\right)^2\notag\\
&<& 2m + \frac{m}{c\binlog n} + \frac{12m}{c} + \frac{4m}{\sqrt{c}}\notag\\
&<& 2m + \frac{17m}{\sqrt{c}}.\label{eqb}
\end{eqnarray}
Thirdly,
\begin{eqnarray*}
\sum_I \hit{I}\length{I} &=& \sum_{\text{base $I$}}\hit{I}\frac{1}{\sqrt{mn}} + \sum_{\substack{\text{regular $I$}\\ \hit{I}=1}}\length{I} \\ && + \sum_{\substack{\text{non-base $I$}\\ \text{$h_I\geq 2$ if regular}}}\hit{I}\length{I}\\
&<& nm\frac{1}{\sqrt{mn}} + n + \left(3n\binlog n +\frac{m}{\sqrt{c}}\right)\left(\frac{2cn\binlog n}{m}\right)\\
&<& \frac{m}{\sqrt{c\binlog n}} + \frac{m}{c\binlog n} + \frac{6m}{c} + \frac{2m}{\sqrt{c}}\\
&<& \frac{10m}{\sqrt{c}}.
\end{eqnarray*}
Thus if $\sqrt{c}\geq 20$, then \eqref{aprequire} is satisfied, and in fact
\begin{equation}\label{eqc}
m-2\sum_I \hit{I}\length{I} > m - \frac{10m}{\sqrt{c}}.
\end{equation}
Substituting estimates \eqref{eqa}, \eqref{eqb}, and \eqref{eqc} into \eqref{apthree} we obtain
\[ \frac{16m}{\sqrt{c}} > \frac{\left(m-\frac{10m}{\sqrt{c}}\right)^2}{2m+\frac{17m}{\sqrt{c}}}, \text{ i.e., } \frac{16}{\sqrt{c}} > \frac{\left(1-\frac{10}{\sqrt{c}}\right)^2}{2+\frac{17}{\sqrt{c}}},\]
which is false if $c$ is sufficiently large, e.g.\ for $\sqrt{c}=100$.

Thus $m<cn\binlog n$, and the proof is finished.
\end{proof}

\section{Final remarks}
\subsection{Infinite dimensions}
In infinite dimensions almost nothing is known.
The most important open question here is the following.
\begin{problem}
Does there exist an infinite equilateral set in any separable infinite-dimensional normed space?
\end{problem}
As observed by Terenzi \cite{MR90c:46017}, it follows from the partition calculus of set theory that a normed space of cardinality at least $2^\mathfrak{c}$ has an infinite equilateral set (in fact an uncountable one), where $c=2^{\aleph_0}$ is the cardinality of the continuum.
The following two easy results are proved in \cite{Swanepoel-thesis}.
\begin{theorem}\label{inf1}
Let $\epsi>0$ and $X$ an infinite dimensional Banach space.
Then $X$ has an equivalent norm with Banach-Mazur distance of at most $2(1+\epsi)$ to the original norm, admitting an infinte equilateral set.
\end{theorem}
\begin{theorem}\label{inf2}
Let $\epsi>0$ and $X$ an infinite dimensional superreflexive Banach space.
Then $X$ has an equivalent norm with Banach-Mazur distance of at most $1+\epsilon$ to the original norm, admitting an infinite equilateral set.
\end{theorem}

\subsection{Generalizations}
An \emph{$\epsi$-almost-equilateral set} in a normed space is a set $S$ satisfying $1-\epsi\leq\norm{\vx-\vy}\leq 1+\epsi$ for all $\vx,\vy\in S$, $\vx\neq \vy$.
Such sets have been studied in \cite{MR86a:46018, MR88h:46033, MR91h:46024, MR93a:46017, MR96j:46006, MR2000b:46013}.
\begin{problem}
Prove that for each $\epsi>0$ there exists $\delta>0$ such that each $n$-dimensional normed space has an $\epsi$-almost-equilateral set of size $(1+\delta)^n$.
\end{problem}
Equilateral sets may be generalized to the notion of a $k$-distance set, i.e., a subset $S$ of a metric space such that $\{d(x,y) : x,y\in S, x\neq y\}$ consists of at most $k$ numbers.
Then the Euclidean case is already non-trivial, and for normed spaces a few results are known: see \cite{MR99k:52028}.
The following conjecture would generalize Petty's theorem:
\begin{conjecture}[\cite{MR99k:52028}]
A $k$-distance set in an $n$-dimensional normed space has size at most $(k+1)^n$.
If equality holds for some $k\geq 1$, then the space must be isometric to $\ell_\infty^n$.
\end{conjecture}
It is known to be true for $n=2$ and for $\ell_\infty^n$, while for general Minkowski spaces there are some weaker estimates; see \cite{MR99k:52028}.

\subsection*{Acknowledgements}
I thank Rafael Villa, Juan Carlos Garc\'{\i}a-V\'azquez, and Jos\'e Antonio Facenda for valuable hints.

\providecommand{\MR}{\relax\ifhmode\unskip\space\fi MR }
\providecommand{\MRhref}[2]{%
  \href{http://www.ams.org/mathscinet-getitem?mr=#1}{#2}
}
\providecommand{\href}[2]{#2}

\end{document}